\documentclass[11pt]{article}

\usepackage[usenames,dvipsnames]{xcolor}

\usepackage[margin={2.5cm,2.2cm}]{geometry} % Define margins

\usepackage{graphicx}
\usepackage{float}
\usepackage{subcaption}
\usepackage{tabularx}

\usepackage{tikz}

\usepackage{enumerate}
\usepackage{amsmath,amsthm,amsfonts,bm}

\usepackage{multirow}
\usepackage{rotating}

\usepackage[backend=bibtex,giveninits=true,sortcites,style=trad-plain]{biblatex}
\usepackage[colorlinks=true]{hyperref}
\hypersetup{
linkcolor=BrickRed
,citecolor=Green
,filecolor=Mulberry
,urlcolor=NavyBlue
,menucolor=BrickRed
,runcolor=Mulberry
,linkbordercolor=BrickRed
,citebordercolor=Green
,filebordercolor=Mulberry
,urlbordercolor=NavyBlue
,menubordercolor=BrickRed
,runbordercolor=Mulberry
}

\theoremstyle{plain}% default
\newtheorem{theorem}{Theorem}[section]
\newtheorem{proposition}[theorem]{Proposition}
\newtheorem{lemma}[theorem]{Lemma}
\newtheorem{corollary}[theorem]{Corollary}
\newtheorem{remark}[theorem]{Remark}

\theoremstyle{definition}

\newtheorem{hypothesis}{Hypothesis}

\usepackage{mathtools}

\definecolor{bananamania}{rgb}{0.98, 0.91, 0.71}
\definecolor{bubbles}{rgb}{0.91, 1.0, 1.0}
\definecolor{lavenderblush}{rgb}{1.0, 0.94, 0.96}
\definecolor{lightred}{RGB}{255,210,202}
\definecolor{lightorange}{RGB}{255,234,199}
\definecolor{lightblue}{HTML}{EEF4F8}
\definecolor{lightyellow}{HTML}{FAF0E1}
\definecolor{lightgreen}{HTML}{F6FFE6}

\newcommand{\R}{\mathbb{R}}

\newcommand{\bv}{\overline{v}}
\newcommand{\bw}{\overline{w}}
\newcommand{\bu}{\overline{u}}

\newcommand{\bbeta}{\bm{\beta}}

\newcommand{\nn}{\mathbf{n}}

\newcommand{\vK}{v_{K}}

\newcommand{\vL}{v_{L}}
\newcommand{\uK}{u_{K}}
\newcommand{\uKs}{u_{K^*}}
\newcommand{\uL}{u_{L}}

\DeclareMathOperator*{\argmin}{arg\,min}

\def\escalar#1#2{\left(#1,#2\right)}

\def\escalarL#1#2{\escalar{#1}{#2}}

\def\escalarML#1#2{\escalar{#1}{#2}_h}

\def\norma#1{\left\|#1\right\|}

\def\salto#1{\left[\!\left[#1\right]\!\right]}
\def\media#1{\left\{\!\!\left\{#1\right\}\!\!\right\}}

\def\T{\mathcal{T}}
\def\E{\mathcal{E}}

\def\N{\mathbb{N}}

\def\X{\mathcal{X}}

\def\Ehi{\mathcal{E}_h^{\text{i}}}
\def\Ehb{\mathcal{E}_h^{\text{b}}}

\def\Pd{\mathbb{P}^{\text{disc}}}
\def\Pc{\mathbb{P}^{\text{cont}}}

\def\Pih{\Pi^h}

\def\aupw#1#2#3{a_h^{\text{upw}}(#1;#2,#3)}

\title{\textbf{
	On a linear DG approximation of chemotaxis models with damping gradient nonlinearities}}

\makeatletter
\def\@fnsymbol#1{\ensuremath{\ifcase#1\or \dagger\or \ddagger\or \mathsection\or * \or\mathparagraph\or *\or **\or \dagger\dagger \or \ddagger\ddagger \else\@ctrerr\fi}}
\makeatother

\author{Daniel Acosta-Soba\thanks{Departamento de Matemáticas, Universidad de Cádiz, Puerto Real, 11510 Cádiz, Spain -- Email: \texttt{daniel.acosta@uca.es} -- Corresponding author}~,
~Alessandro Columbu \thanks{Dipartimento di Matematica e Informatica, Università di Cagliari, Via Ospedale 72, 09124 Cagliari, Italy -- Email: \texttt{alessandro.columbu2@unica.it}}~,
~J. Rafael Rodríguez-Galván\thanks{Departamento de Matemáticas, Universidad de Cádiz, Puerto Real, 11510 Cádiz, Spain -- Email: \texttt{rafael.rodriguez@uca.es}}}
\date{}

\begin{document}

\maketitle

\begin{abstract}
	In this work we present and analyse a novel linear and positivity preserving upwind discontinuous Galerkin (DG) approximation of a class of chemotaxis models with damping gradient nonlinearities. In particular, both a local and a nonlocal model including nonlinear diffusion, chemoattraction, chemorepulsion and logistic growth are considered. Some numerical experiments in the context of chemotactic collapse are presented, whose results are in accordance with the previous analysis of the approximation and show how the blow-up can be prevented by means of the damping gradient term.
\end{abstract}

\paragraph{Keywords:} Upwind. Discontinuous Galerkin. Positivity. Linear scheme. Chemotactic collapse.

\section{Introduction}
\label{sec:chemotaxis_sec}
Since the early work of E.~Keller and L.~Segel \cite{keller_segel_1970,keller_segel_1971} were published, many different alternatives have been proposed trying to extend their ideas. While this early work was intended to model the movement of cells towards high concentration of a chemical signal, other authors have applied this system to modeling population dynamics where the predators are attracted by high concentrations of preys, \cite{kareiva1987swarms}.

These models of population dynamics typically include terms describing population growth or decay in terms of the population size (or density) such as logistic growth terms \cite{mimura1996aggregating,aida2006lower}. In addition, some works suggest that the population of predators can decrease due to accidental deaths as a consequence of their displacements in the chase of their preys, which increase their risk to be decimated by other species that can react to their motion (see \cite{geipel2020predation,visser2007motility,almeda2017behavior}). This accidental deaths have been modeled in terms of the gradient of the predators' population density, \cite{souplet1996finite}.

All the aforementioned effects have been collected in  the model presented in \cite{ishida2023keller}, where the authors are more interested in the mathematical analysis of its solution than in its biological interest. In this sense, the existence and regularity of solutions in chemotaxis models are not clear as, under certain conditions, the so called chemotaxis collapse can occur (see \cite{bellomo2015toward} and the references therein). In this regard, conditions involving the logistic growth and the gradient nonlinearity terms are studied in \cite{ishida2023keller} to ensure the existence of a global, in time, and bounded solution and, therefore, the absence of blow-up.

Regarding the numerical approximation of these models, the literature on positivity preserving and stable numerical approximations of chemotaxis models have experienced a considerable growth during the last years, see, for instance, \cite{acosta-soba_KS_2022,shen_unconditionally_2020,huang_bound_preserving,andreianov_finite_2011,saito_conservative_2007,zhou_finite_2017,chertock_second-order_2008,badia2023bound,guillen2024chemo} and the references therein. Nevertheless, to our best knowledge there is no such numerical scheme for a chemotaxis model involving a gradient damping term. This term presents an additional remarkable difficulty in the context of chemotaxis collapse, where very steep gradients of the cell density appear, possibly leading to spurious oscillations if it is not treated adequately.

This work aims to be complementary to the paper \cite{li2025dissipative}, where the authors introduce an even more general version of the model in \cite{ishida2023keller}, including a damping gradient term. In this previous paper, some numerical experiments in the context of blowing-up solutions are carried out using the linear, positivity preserving numerical approximation that we are going to develop in the following sections.

The structure of this document is outlined as follows. In Section~\ref{sec:chemotaxis} we present the two different extensions of the model in \cite{ishida2023keller} that are studied in \cite{li2025dissipative}. On the one hand, a local model that additionally includes possibly nonlinear diffusion, chemoattraction and chemorepulsion is described in \eqref{problem:chemotaxis_local}. On the other hand, a similar model including nonlocal effects of the cell density into the chemosensitivities can be considered, see \eqref{problem:chemotaxis_nonlocal}. Some results regarding the local and global in time existence and regularity of the solutions of both the local and nonlocal models under certain constraints on the parameters that have been developed in \cite{li2025dissipative} are stated in Section~\ref{sec:results} without proofs.

Then, in Section~\ref{sec:esquema_completamente_discreto} we present and analyze the properties of the linear, positivity preserving, discrete approximation of the models used in the numerical experiments shown in \cite{li2025dissipative}. This numerical approximation is based on the upwind discontinuous Galerkin ideas previously introduced in \cite{acosta-soba_upwind_2022,acosta-soba_KS_2022}. Since the numerical test in \cite{li2025dissipative} were focused on the parabolic-elliptic ($\tau=0$) local model \eqref{problem:chemotaxis_local}, in Section~\ref{sec:numer-experiments} we show some complementary numerical experiments for the fully parabolic ($\tau=1$) local model~\eqref{problem:chemotaxis_local} and the nonlocal model~\eqref{problem:chemotaxis_nonlocal}. Similarly to those in \cite{li2025dissipative}, these results, which are in accordance with the previous theoretical analysis, show that the suppression of some of the analytical conditions on the parameters (and in particular those associated to the
gradient nonlinearities \eqref{condgamma}) that ensure global in time existence of the solution may lead to only local in time solutions that exhibit finite time blow-up.

\section{Chemotaxis models with damping gradient nonlinearities}
\label{sec:chemotaxis}

Let $\Omega\subset\R^d$ with $d\in\N$ be a regular-enough bounded domain, e.g. polygonal (in case $d=2$), with $\nn$ its unit normal vector and $T>0$.

The classical Keller-Segel model consists of finding two real valued functions: $u=u(x,t)$, the cell density, and $v=v(x,t)$, the chemoattractant signal, defined in $\Omega\times [0,T]$ satisfying
\begin{subequations}
	\label{problem:KS}
	\begin{align}
	\label{eq:KS_u}
	u_t&=k_0\Delta u-k_1\nabla\cdot(u\nabla v),&\text{in }\Omega\times (0,T),\\
	\label{eq:KS_v}
	\tau v_t&=k_2\Delta v-k_3v+k_4u,&\text{in }\Omega\times (0,T),\\
	0&=\nabla u\cdot \mathbf{n}=\nabla v\cdot\mathbf{n} 
	, &\text{on }\partial\Omega\times (0,T),\\
	\label{eq:CI.cahn-hilliard+adveccion}
	u(\cdot,0)&=u_0,\quad v(\cdot,0)=v_0\text{ if }\tau>0,&\text{in }\Omega.%\\
	\end{align}
\end{subequations}
 
In this section, first, we consider a more general version of the classical model \eqref{problem:KS} which introduces nonlinear diffusion, attraction and repulsion, logistic growth and damping gradient nonlinearities. In particular, this \textbf{local model}, which has been first introduced in \cite{li2025dissipative}, 
extends the one in \cite{ishida2023keller} and
consists of finding three real valued functions: $u=u(x,t)$, the cell density, $v=v(x,t)$, the chemoattractant signal, and $w=w(x,t)$, the chemorepulsive signal, defined
in $\Omega\times [0,T]$ satisfying
\begin{subequations}
	\label{problem:chemotaxis_local}
	\begin{align}
		\label{eq:chemotaxis_local_u}
		u_t=& \nabla \cdot \left((u+1)^{n_1-1}\nabla u\right) -\chi \nabla \cdot \left(u(u+1)^{n_2-1}\nabla v\right) \nonumber\\&+\xi \nabla \cdot \left(u(u+1)^{n_3-1}\nabla w\right)+\lambda u^\rho -\mu u^k - c\vert\nabla u\vert^\gamma & \text{in } \Omega \times (0,T),\\
		\label{eq:chemotaxis_local_v}
		\tau v_t=& \Delta v - av  + f_1(u) & \text{in } \Omega \times (0,T),\\
		\label{eq:chemotaxis_local_w}
		\tau w_t=& \Delta w -  dw + f_2(u) & \text{in } \Omega \times (0,T),\\
		0=&\nabla u\cdot\nn=\nabla v\cdot\nn=\nabla w\cdot\nn & \text{on } \partial \Omega \times (0,T),\\
		u(\cdot,0)=&u_0(x), \tau v(\cdot,0)=\tau v_0(x), \tau w(\cdot,0)=\tau w_0(x)  & \text{in } \Omega.
	\end{align}
\end{subequations}
Moreover, we consider a \textbf{nonlocal model} where $v=v(x,t)$ and $w=w(x,t)$ represent the chemoattractant and the chemerepulsive deviations, respectively, satisfying
\begin{subequations}
	\label{problem:chemotaxis_nonlocal}
	\begin{align}
		\label{eq:chemotaxis_nonlocal_u}
		u_t=& \nabla \cdot \left((u+1)^{n_1-1}\nabla u\right) -\chi \nabla \cdot \left(u(u+1)^{n_2-1}\nabla v\right) \nonumber\\&+ \xi \nabla \cdot \left(u(u+1)^{n_3-1}\nabla w\right)+\lambda u^\rho -\mu u^k - c\vert\nabla u\vert^\gamma & \text{in } \Omega \times (0,T),\\
		\label{eq:chemotaxis_nonlocal_v}
		0=& \Delta v - \frac{1}{\vert\Omega\vert}\int_\Omega f_1(u)  + f_1(u) & \text{in } \Omega \times (0,T),\\
		\label{eq:chemotaxis_nonlocal_w}
		0=& \Delta w - \frac{1}{\vert\Omega\vert}\int_\Omega f_2(u)  + f_2(u) & \text{in } \Omega \times (0,T),\\
		0=&\nabla u\cdot\nn=\nabla v\cdot\nn=\nabla w\cdot\nn & \text{on } \partial \Omega \times (0,T),\\
		u(\cdot,0)=&u_0(x)  & \text{in } \Omega,\\
		\label{eq:chemotaxis_nonlocal_const}
		0=&\int_\Omega v(x,t)dx=\int_\Omega w(x,t)dx & \text{in } (0,T).
	\end{align}
\end{subequations}
Here, $\tau\in\{0,1\}$, $\chi,\xi,\lambda,\mu,c\ge 0$, $n_1,n_2,n_3\in\R$ and $\rho, k, \gamma\ge 1$.

Notice that the local model \eqref{problem:chemotaxis_local} is a direct generalization of the classical Keller--Segel model \eqref{problem:KS} where $v=v(x,t)$ represents the chemoattractant density and a new variable concerning the chemorepulsion density, $w=w(x,t)$, has been considered. In particular, equations \eqref{eq:chemotaxis_local_u} or \eqref{eq:chemotaxis_nonlocal_u} coincide with the equation \eqref{eq:KS_u} for $n_1=n_2=1$ and $\xi=\lambda
=\mu=c=0$. In this sense, equations \eqref{eq:chemotaxis_local_u} or \eqref{eq:chemotaxis_nonlocal_u} include possibly nonlinear diffusion and chemoattraction, in the case $n_1,n_2\neq 1$, and chemorepulsion if $\xi\neq 0$. In addition, a logistic growth term $\lambda u^\rho-\mu u^k$
and a damping gradient term $- c\vert\nabla u\vert^\gamma$
 have been added to equations \eqref{eq:chemotaxis_local_u} and \eqref{eq:chemotaxis_nonlocal_u}.

In the considered nonlocal model \eqref{problem:chemotaxis_nonlocal}, the unknown $v=v(x,t)$ represents  the chemoattractant deviation and not its density as in the local model \eqref{problem:chemotaxis_local} or in the classical Keller-Segel model \eqref{problem:KS}. In particular, since the deviation measures the difference between the observed value of a specific variable and its mean, we have that $v$ changes sign, contrarily  to what occurs to the cell and signal densities (which are nonnegative). Subsequently, the mean of $v$ is zero, exactly as imposed  in the equation \eqref{eq:chemotaxis_nonlocal_const}. Naturally, the chemical repellent $w$ behaves similarly. In this sense, we do not make use of a different symbolism for the deviations, since from the context it is clear to which quantity we are referring to. 

  \begin{remark}
    Notice that we assume $\rho,k,\gamma\ge 1$ so that the terms involving these parameters in \eqref{eq:chemotaxis_local_u} or in \eqref{eq:chemotaxis_nonlocal_u} are differentiable when $u=0$ or $\vert\nabla u\vert=0$. 
  \end{remark}

In this section, given 
$\delta\in(0,1)$, we denote by $C^{k+\delta}(\overline\Omega)$ the set of functions that are $k$-times, with $k\in\N\cup\{0\}$, continuously differentiable and whose $k$-th derivative satisfies a Hölder condition with index $\delta$ (see e.g. \cite{evans2010partial}).  Moreover, $C^{k+\delta,m+\eta}(\overline\Omega\times[0,T])$ is defined for certain $\delta,\eta\in(0,1)$, $k,m\in\N\cup\{0\}$, by functions $u(x,t)$ which satisfy
$u(\cdot, t)\in C^{k+\delta}(\overline\Omega)$ and $u(x,\cdot)\in C^{m+\eta}([0,T])$ for every $t\in[0,T]$ and $x\in\overline\Omega$.

Following \cite{li2025dissipative} the following assumptions are made so that the results in Section~\ref{sec:results} hold.

\begin{hypothesis}
	\label{hyp:reglocal}
Here, we assume that $\Omega$ is an open bounded domain of class $C^{2+\delta}$ for some $\delta\in(0,1)$ (see \cite{evans2010partial} for details), the sources $f_1, f_2$ are such that $f_1,f_2: [0,\infty) \rightarrow \R_+$,
 with  $f_1,f_2\in C^1([0,\infty))$ and the initial data $u_0=u_0(x), \tau v_0=\tau v_0(x)$ and $\tau w_0=\tau w_0(x)$ satisfy $u_0, v_0, w_0: \overline{\Omega}  \rightarrow \R_+$,
  with 
$$u_0,  v_0,   w_0 \in 
\{\psi \in C^{2+\delta}(\overline{\Omega}): \nabla\psi\cdot\nn=0 \text{ on } \partial \Omega\}.$$
\end{hypothesis}
\begin{hypothesis}
	\label{hyp:disf}
  We assume that for proper $\alpha,\beta>0$,
\begin{equation}
    f_1(s)\leq k_1 (s+1)^\alpha \quad \text{and} \quad  k_2(s+1)^\beta \leq f_2(s)\leq k_3 (s+1)^\beta,
\end{equation}
where $k_1,k_2,k_3> 0$ are positive constants.
\end{hypothesis}

\begin{remark}
  \label{rmk:regularity_f}
  The bounds on $f_1$ and $f_2$ in Hypothesis~\ref{hyp:disf} are required for technical reasons in the proofs of Theorems~\ref{thm:theoremlocal} and \ref{thm:theoremnonlocal} (see \cite{li2025dissipative}).
   In fact, these bounds can be relaxed to the following:
  \begin{equation*}
    f_1(s)\leq k_1 (s+\eta)^\alpha \quad \text{and} \quad  k_2(s+\eta)^\beta \leq f_2(s)\leq k_3 (s+\eta)^\beta,
\end{equation*}
for $\eta>0$. Notice that we avoid the case $\eta= 0$ to move away from a possible singularity at $s=0$ of the functions $f_1(s)=s^\alpha$ and $f_2(s)=s^\beta$ as they have to be in $C^1([0,\infty))$.
\end{remark}

\section{Existence and regularity of solutions}
\label{sec:results}
Now, we state the conditions that ensure us the global boundedness of the solutions. In particular, we want to bound solutions which blow up by adding the gradient term $-c\vert\nabla u\vert^\gamma$ to the logistic and, therefore, we do not focus on the bounding effect that the diffusive term $\nabla \cdot \left((u+1)^{n_1-1}\nabla u\right)$ has. In order to do that, we establish the following relation involving the parameters in the local and the nonlocal models, \eqref{problem:chemotaxis_local} and \eqref{problem:chemotaxis_nonlocal}, respectively: for any $d\in\N$, $\tau\in\{0,1\}$, $n_1,n_2,n_3\in\R$ and 
$\alpha,\beta,\chi,\xi,\lambda,\mu,c>0$, $1\leq\rho< k$ and $\gamma\in[1,2]$ we fix 
\begin{equation} \label{condgamma}  
\max{\left\{1,\frac{d}{d+1}(n_2+\alpha),\tau\frac{d}{d+1}(n_3+\beta)\right\}}<\gamma\leq2. 
\end{equation}

\begin{remark}
  Here, in condition \eqref{condgamma}, we have assumed that $1\leq\rho< k$ and $\gamma\in[1,2]$. These two assumptions are required for the statements in Theorems~\ref{thm:theoremlocal} and \ref{thm:theoremnonlocal}.

  On the one hand, $1\leq\rho< k$ implies that $\lambda u^\rho -\mu u^k$ in \eqref{eq:chemotaxis_local_u} or \eqref{eq:chemotaxis_nonlocal_u} is a logistic term so that the decay becomes stronger than the production when $u$ increases. This assumption, along with the negative sign of the term $-c\vert\nabla u\vert^\gamma$, allows us to have a bound on the mass of $u$ as in Proposition~\ref{prop:mass_bound}.

  On the other hand, we assume $\gamma\in[1,2]$ for technical reasons concerning the results in Theorems~\ref{thm:theoremlocal} and \ref{thm:theoremnonlocal} when $c>0$. Intuitively, the idea of this assumption is that, if $\gamma >2$, the solution might not be global in time due to $\norma{\nabla u}_{L^\infty(\Omega)^d}$ not being bounded in finite time even though $\norma{u}_{L^\infty(\Omega)}$ is bounded. However, if $\gamma\in[1,2]$, the previous case cannot occur and if the solution is only local in time, i.e, it exists up to some finite time $T$, then necessarily $\norma{u}_{L^\infty(\Omega)}$ blows up at $T$.
  
  More details can be consulted in \cite{li2025dissipative,ishida2023keller} and the references therein.
\end{remark}

The succeeding theorems are stated without proofs.
The full discussion and proofs of these 
results can be found in \cite{li2025dissipative}.
\begin{proposition}
  \label{prop:mass_bound}
	Let $\mu>0$ and $1\le\rho<k$. The mass of the cell density $u$ in \eqref{problem:chemotaxis_local} or \eqref{problem:chemotaxis_nonlocal} is bounded as follows:
	\begin{equation}
		\label{mass_bound}
	\int_\Omega u \le\max\left\{\int_\Omega u_0,\left(\frac{\lambda}{\mu}\vert\Omega\vert^{k-\rho}\right)^{\frac{1}{k-\rho}}\right\}, \quad\forall t\in(0,T),
	\end{equation}
  where $T$ is the upper bound of the time interval where the function $u$ is defined, with $T\le\infty$.
\end{proposition}
\begin{theorem}[Local classical solution of \eqref{problem:chemotaxis_local} and \eqref{problem:chemotaxis_nonlocal}]\label{thm:theoremExistence}
For $\tau\in \{0,1\}$ and $\delta\in(0,1)$, let $\Omega$, $f_1,f_2, u_0, \tau v_0,$ $\tau w_0$ comply with Hypothesis~\ref{hyp:reglocal}.  Additionally, let $\chi,\xi,\lambda,\mu>0$, $n_1,n_2,n_3\in \R$, $k,\rho,\gamma \geq 1$.  Then there exist $T>0$ and a unique triple of functions $(u,v,w)$, with 
\begin{equation*}
(u,v,w)\in C^{2+\delta,1+\frac{\delta}{2}}( \overline{\Omega} \times [0, T])\times C^{2+\delta,\tau+\frac{\delta}{2}}( \overline{\Omega} \times [0, T])\times  C^{2+\delta,\tau+\frac{\delta}{2}}( \overline{\Omega} \times [0, T]),
\end{equation*}
solving problems  \eqref{problem:chemotaxis_local} and \eqref{problem:chemotaxis_nonlocal}, where $u\ge 0$ in $\overline\Omega \times [0,T]$ and, in the case \eqref{problem:chemotaxis_local}, $v,w\ge 0$ in $\overline\Omega \times [0,T]$.
\end{theorem}
\begin{theorem}[Global and bounded classical solution of \eqref{problem:chemotaxis_local}]\label{thm:theoremlocal}
For $\tau\in \{0,1\}$, $\delta\in(0,1)$ and $\alpha,\beta>0$, let $\Omega$, $f_1,f_2, u_0$ comply with Hypotheses~\ref{hyp:reglocal} and \ref{hyp:disf}.  Additionally, let $\chi,\xi,\lambda,\mu,c>0$, $n_1,n_2,n_3\in\R$, $1\leq\rho< k$ and condition \eqref{condgamma} hold true. Then problem \eqref{problem:chemotaxis_local} admits a unique solution
\begin{equation*}
(u,v,w)\in C^{2+\delta,1+\frac{\delta}{2}}( \overline{\Omega} \times [0, \infty))\times C^{2+\delta,\tau+\frac{\delta}{2}}( \overline{\Omega} \times [0, \infty))\times  C^{2+\delta,\tau+\frac{\delta}{2}}( \overline{\Omega} \times [0, \infty)) 
\end{equation*}
such that $u,v,w\geq 0$ in $\overline\Omega \times [0,\infty)$ and $u,v,w \in L^\infty(\Omega \times (0,\infty)).$
\end{theorem}
\begin{theorem}[Global and bounded classical solution of \eqref{problem:chemotaxis_nonlocal}]\label{thm:theoremnonlocal}
For $\delta\in(0,1)$ and $\alpha,\beta>0$,let $\Omega$, $f_1,f_2, u_0$ comply with Hypotheses~\ref{hyp:reglocal} and \ref{hyp:disf}.
 Additionally, let $\alpha,\beta,\chi,\xi,\lambda,\mu,c>0$, $n_1,n_2,n_3\in\R$, $1\leq\rho< k$ and condition \eqref{condgamma} hold true.  Then problem \eqref{problem:chemotaxis_nonlocal} admits a unique solution 
 \begin{equation*}
(u,v,w)\in C^{2+\delta,1+\frac{\delta}{2}}( \overline{\Omega} \times [0, \infty))\times C^{2+\delta,\frac{\delta}{2}}( \overline{\Omega} \times [0, \infty))\times  C^{2+\delta,\frac{\delta}{2}}( \overline{\Omega} \times [0, \infty)) 
\end{equation*}
such that $u\geq 0$ in $\overline\Omega \times [0,\infty)$ and $u,v,w \in L^\infty(\Omega \times (0,\infty)).$
\end{theorem}

\begin{remark}
    Notice that the values in the condition \eqref{condgamma} are not critical. In fact, the influence of $n_1$ in this condition has not been taken into consideration in this work. In this sense, it has been studied in previous works (see, for instance, \cite{ColumbuFrassuViglialoro2023,li2023combining,columbu2023uniform}) how the diffusive term $\nabla\cdot\left((u+1)^{n_1-1}\nabla u\right)$ in \eqref{eq:chemotaxis_local_u} or in \eqref{eq:chemotaxis_nonlocal_u} affects the global in time existence of solution of similar models. 
\end{remark}

\section{Numerical approximation}
\label{sec:esquema_completamente_discreto}

We propose a linear, computationally efficient numerical scheme to approximate the local and the nonlocal models, \eqref{problem:chemotaxis_local} and \eqref{problem:chemotaxis_nonlocal}, that provides a physically meaningful approximation of every variable.

\subsection{Notation}

 We consider a finite element shape-regular triangular mesh $\T_h=\{K\}_{K\in \T_h}$ in the sense of Ciarlet, \cite{ciarlet2002finite}, of size $h$ over $\Omega$. We denote by $\E_h$ the set of the edges of $\T_h$ (faces if $d=3$) with $\Ehi$ the set of the \textit{interior edges} and $\Ehb$ the \textit{boundary edges}, thus $\E_h=\Ehi\cup\Ehb$.

 Now, we fix the following orientation over the mesh $\T_h$:
 \begin{itemize}
	\item For any interior edge $e\in\Ehi$ we set the associated unit normal vector $\nn_e$. In this sense, when refering to edge $e\in\Ehi$ we will denote by $K_e$ and $L_e$ the elements of $\T_h$ with $e=\partial K_e\cap\partial L_e$ and so that $\nn_e$ is exterior to $K_e$ pointing to $L_e$.
	
	If there is no ambiguity, to abbreviate the notation  we will denote the previous elements $K_e$ and $L_e$ simply by $K$ and $L$, respectively, with the assumption that their naming is always with respect to the edge $e\in\Ehi$ and it may vary if we consider a different edge of $\Ehi$.
	\item For any boundary edge $e\in\Ehb$, the unit normal vector $\nn_e$ points outwards of the domain $\Omega$.
 \end{itemize}

 Therefore, we can define the \textit{average} $\media{\cdot}$ and the \textit{jump} $\salto{\cdot}$ of a function $v$ on an edge $e\in\E_h$ as follows:
\begin{equation*}
		\media{v}\coloneqq
		\begin{cases}
			\dfrac{\vK+\vL}{2}&\text{if } e\in\Ehi,\,  e=K\cap L,\\
			\vK&\text{if }e\in\Ehb,\, e\subset K,
		\end{cases}
		\salto{v}\coloneqq
		\begin{cases}
			\vK-\vL&\text{if } e\in\Ehi,\, e=K\cap L,\\
			\vK&\text{if }e\in\Ehb,\, e\subset K.
		\end{cases}
\end{equation*}

We denote by $\Pc_k(\T_h)$ and $\Pd_k(\T_h)$ the spaces of continuous and discontinuous finite element functions, respectively, which are polynomials of degree $k\ge0$ when restricted to the elements $K$ of $\T_h$. See, for instance, \cite{ern_theory_2010,di_pietro_mathematical_2012}.

Moreover, we take a partition $0=t_0<t_1<\cdots<t_N=T$ of the time domain $[0,T]$
with $\Delta t^{m+1}=t_{m+1}-t_m$ the time step. For any function $v$ depending on time, we denote 
$v^{m+1}\simeq v(t_{m+1})$ and the discrete backward Euler time derivative operator $v_t(t_{m+1})\simeq\delta_t v^{m+1}:=(v^{m+1}-v^m)/\Delta t^{m+1}$.

In addition, we set the following notation for the positive and negative parts of a function $v$:
	$$
	v_\oplus\coloneqq\frac{|v|+v}{2}=\max\{v,0\},
	\quad
	v_\ominus\coloneqq\frac{|v|-v}{2}=-\min\{v,0\},
	\quad
	v=v_\oplus - v_\ominus.
	$$

Finally, we define the projection operators $\Pi_1\colon L^1(\Omega)\longrightarrow \Pc_1(\T_h)$ and 
$\Pih_1\colon L^1(\Omega)\longrightarrow \Pc_1(\T_h)$ 
as follows:
\begin{align}
  \label{eq:esquema_DG_Pi1}
  \escalarL{\Pi_1 g}{\overline{v}}
  &=
  \escalarL{g}{\overline{v}},&\forall\,\overline{v}\in \Pc_1(\T_h),
  \\
  \label{eq:esquema_DG_Pih1}
  \escalarML{\Pih_1 g}{\overline{v}}
&=\escalarL{g}{\overline{v}},&\forall\,\overline{v}\in \Pc_1(\T_h),
\end{align}
where $\escalarL{\cdot}{\cdot}$ denotes the usual scalar product in $L^2(\Omega)$ and $\escalarML{\cdot}{\cdot}$ denotes the mass-lumping scalar product in $\Pc_1(\T_h)$ resulting from using the trapezoidal rule to approximate the scalar product in $L^2(\Omega)$ (see, for instance, \cite{quarteroni2008numerical}).

\subsection{Fully discrete scheme}

In each time step $m+1$, we decouple the equations for the chemical signals \eqref{eq:chemotaxis_local_v}--\eqref{eq:chemotaxis_local_w} or \eqref{eq:chemotaxis_nonlocal_v}--\eqref{eq:chemotaxis_nonlocal_w} from the equations for the cell density, \eqref{eq:chemotaxis_local_u} or \eqref{eq:chemotaxis_nonlocal_u}, respectively. In this way, we first compute
a $\Pc_1(T_h)$-approximation of $v$ and $w$ treating $u$
 explicitly in \eqref{eq:chemotaxis_local_v}--\eqref{eq:chemotaxis_local_w} or in \eqref{eq:chemotaxis_nonlocal_v}--\eqref{eq:chemotaxis_nonlocal_w}, where $u^m\ge0$ as shown in Section~\ref{sec:properties_scheme}.

 In fact, we can obtain the following linear finite element scheme
 \begin{align*}
  \tau\escalarML{\delta_tv^{m+1}}{\bv}&=-\int_\Omega \nabla v^{m+1}\cdot\nabla\bv-a\escalarML{v^{m+1}}{\bv}+\escalarL{f_1(u^m)}{\bv},\\
 	\tau\escalarML{\delta_tw^{m+1}}{\bw}&=-\int_\Omega\nabla w^{m+1}\cdot\nabla\bw-d\escalarML{w^{m+1}}{\bw}+\escalarL{f_2(u^m)}{\bw},
 \end{align*}
 for all $\bv,\bw\in\Pc_1(\T_h)$, with a unique nonnegative solution that approximates $v$ and $w$ in \eqref{problem:chemotaxis_local} in a mesh of non-obtuse triangles using the well known mass-lumping technique, see, for instance, \cite{acosta-soba_KS_2022}.

On the other hand, we can compute an approximation of $v$ and $w$ in \eqref{problem:chemotaxis_nonlocal} using standard finite elements as follows,
\begin{align*}
  \int_\Omega \nabla v^{m+1}\cdot\nabla\bv&=-\frac{1}{\vert\Omega\vert}\escalarL{\int_\Omega f_1(u^{m})}{\bv}+\escalarL{f_1(u^m)}{\bv},\\
 	\int_\Omega\nabla w^{m+1}\cdot\nabla\bw&=-\frac{1}{\vert\Omega\vert}\escalarL{\int_\Omega f_2(u^{m})}{\bw}+\escalarL{f_2(u^m)}{\bw},
 \end{align*}
 for all $\bv,\bw\in\Pc_1(\T_h)$, where the constraint \eqref{eq:chemotaxis_nonlocal_const} can be enforced by postprocessing.
 
 Thus, hereafter, we will only focus on the approximation of the cell density $u$.

Next, with the approximations that we have obtained for $v$ and $w$, with $v^{m+1},w^{m+1}\ge 0$ in the case \eqref{problem:chemotaxis_local} or $\int_\Omega v^{m+1}=\int_\Omega w^{m+1}=0$ in the case \eqref{problem:chemotaxis_nonlocal}, we compute an approximation of $u$ in the current time step $m+1$. Our aim is to develop a linear positive approximation of $u$
using the ideas previously introduced in \cite{acosta-soba_KS_2022} where the equation of the cell density is rewritten as gradient flux using the chemical potential for the particular case of the classical Keller-Segel model. However, although this is not possible for this chemotaxis variant, we proceed as follows: formally,
$$
\nabla\log(u)=\frac{1}{u}\nabla u,
$$
so that we can rewrite, for $u\ge 0$,
$$
\nabla\cdot((u+1)^{n_1-1}\nabla u)=\nabla\cdot(u(u+1)^{n_1-1}\nabla\log(u)) \quad\text{and}\quad \vert\nabla u\vert=u\vert \nabla\log(u)\vert.
$$
Hence, regularizing the term $\log(u)$ by adding a small parameter $\varepsilon>0$ and projecting it into the space $\Pc_1(\T_h)$ by means of the projection operator $\Pi_1$ defined in \eqref{eq:esquema_DG_Pi1}, we propose the following implicit-explicit linear discrete scheme:

Given $u^m\in\Pd_0(\T_h)$ with $u^m\ge 0$ and $v^{m+1},w^{m+1}\in\Pc_1(\T_h)$, find $u^{m+1}\in \Pd_0(\T_h)$ with $u^{m+1}\ge 0$ solving the problem
\begin{align}
\label{scheme_DG_upw_chemo}
\escalarL{\delta_tu^{m+1}}{\bu}&+\aupw{-(u^m+1)^{n_1-1}\nabla\Pi_1(\log(u^m+\varepsilon))}{u^{m+1}}{\bu}\nonumber\\&+\aupw{(u^m+1)^{n_2-1}\nabla v^{m+1}}{u^{m+1}}{\bu}\nonumber\\&+\aupw{-(u^m+1)^{n_3-1}\nabla w^{m+1}}{u^{m+1}}{\bu}\nonumber\\&-\lambda \escalarL{(u^{m})^\rho}{\bu} +\mu \escalarL{u^{m+1}(u^m)^{k-1}}{\bu}\nonumber\\&+ c \escalarL{u^{m+1} (u^m)^{\gamma-1}\vert\nabla\Pi_1\log(u^m+\varepsilon)\vert^\gamma}{\bu}=0,
\end{align}
for all $\bu\in \Pd_0(\T_h)$, where
\begin{equation}
	\label{def:aupw}
	\aupw{\bbeta}{u}{\bu}\coloneqq \sum_{e\in\E_h^i,e=K\cap L}\int_e\left(\left(\media{\bbeta}\cdot\nn_e\right)_\oplus\uK-(\media{\bbeta}\cdot\nn_e)_\ominus\uL\right)\salto{\bu}.
\end{equation}
Notice that the bilinear form $\aupw{\bbeta}{u}{\bu}$ contains no derivative terms because they vanish due to the selection of piecewise constant functions, see~\cite{acosta-soba_KS_2022}.

The numerical scheme \eqref{scheme_DG_upw_chemo} is only one possible linear approach to obtain a positive approximation of the cell density as shown in the next Section~\ref{sec:properties_scheme}. However, other smart modifications of the proposed approach might be considered to preserve other properties of the continuous models such as the mass bound \eqref{mass_bound}.

In Section~\ref{sec:properties_scheme} we will provide a way of computing the solution of \eqref{scheme_DG_upw_chemo} enforcing $u^{m+1}\ge 0$.

\subsubsection{Properties of the scheme}
\label{sec:properties_scheme}
Consider the following auxiliary nonlinear truncated scheme:
Given $u^m\in\Pd_0(\T_h)$ with $u^m\ge 0$ and $v^{m+1},w^{n+1}\in\Pc_1(\T_h)$, find $u^{m+1}\in \Pd_0(\T_h)$ solving the problem
\begin{align}
\label{scheme_DG_upw_chemo_truncated}
\escalarL{\delta_tu^{m+1}}{\bu}&+\aupw{-(u^m+1)^{n_1-1}\nabla\Pi_1(\log(u^m+\varepsilon))}{u^{m+1}_\oplus}{\bu}\nonumber\\&+\aupw{(u^m+1)^{n_2-1}\nabla v^{m+1}}{u^{m+1}_\oplus}{\bu}\nonumber\\&+\aupw{-(u^m+1)^{n_3-1}\nabla w^{m+1}}{u^{m+1}_\oplus}{\bu}\nonumber\\&-\lambda \escalarL{(u^{m})^\rho}{\bu} +\mu \escalarL{u^{m+1}_\oplus(u^m)^{k-1}}{\bu}\nonumber\\&+ c \escalarL{u^{m+1}_\oplus(u^m)^{\gamma-1}\vert\nabla\Pi_1\log(u^m+\varepsilon)\vert^\gamma}{\bu}=0,
\end{align}
for all $\bu\in \Pd_0(\T_h)$.

\begin{lemma}[Local mass bounds]
	\label{prop:discrete_mass_bounds}
	Assume that there is a solution $u^{m+1}$ of \eqref{scheme_DG_upw_chemo_truncated}, then it satisfies
	\begin{equation}
		\label{mass_bound_discrete}
		\int_\Omega u^{m+1}\le \int_\Omega u^m+\Delta t^{m+1} \lambda\int_\Omega (u^m)^\rho.
	\end{equation}
\end{lemma}
\begin{proof}
	Just test \eqref{scheme_DG_upw_chemo_truncated} by $\overline{u}=1$.
\end{proof}

Notice that the bound of the mass of the discrete solution \eqref{mass_bound_discrete} is consistent with the bound on the mass of the continuous solution \eqref{mass_bound}. 

\begin{theorem}[DG scheme \eqref{scheme_DG_upw_chemo_truncated} preserves positivity]
	\label{thm:positivity_DG_chemo_truncated}
	
	If we assume that $u^m\ge 0$ then any solution of \eqref{scheme_DG_upw_chemo_truncated} satisfies that $u^{m+1}\ge 0$ in $\Omega$.
\end{theorem}
\begin{proof}
	To prove that if $u^m\ge 0$ then $u^{m+1}\ge 0$ just take the following test function in \eqref{scheme_DG_upw_chemo_truncated}:
	$$
	\bu=
	\begin{cases}
		(u_{K^*}^{m+1})_\ominus &\text{in }K^*\in\T_h,\\
		0 & \text{otherwise},
	\end{cases}
	$$
	where $K^*=\argmin_{K\in\T_h}u^{m+1}_K$.

	Therefore, we arrive at 
	\begin{align*}
		\escalarL{\delta_tu^{m+1}}{\bu}&+\aupw{-(u^m+1)^{n_1-1}\nabla\Pi_1(\log(u^m+\varepsilon))}{u^{m+1}_\oplus}{\bu}\\&+\aupw{(u^m+1)^{n_2-1}\nabla v^{m+1}}{u^{m+1}_\oplus}{\bu}\\&+\aupw{-(u^m+1)^{n_3-1}\nabla w^{m+1}}{u^{m+1}_\oplus}{\bu}\\&\ge\lambda\escalarL{(u^m)^\rho}{\bu}\ge0,
	\end{align*}

    Now, since the positive part is a non-decreasing function and for all $L\in\T_h$ we have that $\uL^{m+1}\ge \uKs^{m+1}$ implies $(\uL^{m+1})_{\oplus}\ge (\uKs^{m+1})_{\oplus}$. Consequently, we deduce that
	\begin{align*}
		\aupw{\bbeta}{u^{m+1}_{\oplus}}{\bu}&=\sum_{e\in\Ehi, e=K^*\cap L}\int_e\left((\bbeta\cdot\nn_{e})_{\oplus}(\uKs^{m+1})_{\oplus}-(\bbeta\cdot\nn_{e})_{\ominus}(\uL^{m+1})_{\oplus}\right)(\uKs^{m+1})_{\ominus}\\
		&\le\sum_{e\in\Ehi, e= K^*\cap L}\int_e\left((\bbeta\cdot\nn_{e})_{\oplus}(\uKs^{m+1})_{\oplus}-(\bbeta\cdot\nn_{e})_{\ominus}(\uKs^{m+1})_{\oplus}\right)(\uKs^{m+1})_{\ominus}\\
		&=\sum_{e\in\Ehi, e=K^* \cap L}\int_e(\bbeta\cdot\nn_{e})(\uKs^{m+1})_{\oplus}(\uKs^{m+1})_{\ominus}=0,
	\end{align*}
    for any $\bbeta\in\Pd_1(\T_h)$.
    
    Hence,
    $$
	0\le \escalarL{\delta_tu^{m+1}}{\bu}=\vert K^*\vert \left(\delta_t \uKs^{m+1}\right)(\uKs^{m+1})_{\ominus} =
	-\frac{\vert K^*\vert }{\Delta t^{m+1}}\left((\uKs^{m+1})_{\ominus}^2+\uKs^m(\uKs^{m+1})_{\ominus}\right) \le 0,
	$$
	which implies $(\uKs^{m+1})_{\ominus}=0$.
    
    Thanks to the choice of $K^*$, we can conclude that $u^{m+1}\ge0$.
\end{proof}

To prove the next result we will use the Leray-Schauder fixed point theorem:
  \begin{quote}
    \begin{theorem}[Leray-Schauder fixed point theorem]
        \label{thm:Leray-Schauder}
        Let $\X$ be a Banach space and let $T\colon\X\longrightarrow\X$ be a continuous and compact operator. If the set $$\{x\in\X\colon x=\alpha \,T(x)\quad\text{for some } 0\le\alpha\le1\}$$ is bounded (uniformly with respect to $\alpha$), then $T$ has  at least one fixed point.
    \end{theorem}
\end{quote}

\begin{theorem}
	\label{prop:existence_solution_truncated}
	There is at least one solution of \eqref{scheme_DG_upw_chemo_truncated}.
\end{theorem}
\begin{proof}
	
	Given $u^m\in\Pd_0(\T_h)$ with $u^m\ge 0$ and and $v^{m+1},w^{m+1}\in\Pc_1(\T_h)$, we define the map
	$$T\colon \Pd_0(\T_h)\longrightarrow\Pd_0(\T_h)$$ such that $T(\widehat{u})=u\in\Pd_0(\T_h)$ is the unique solution of the linear problem:

	\begin{align}
	\label{eq:existence_solution_truncated}
	\frac{1}{\Delta t^{m+1}}\escalarL{u-u^m}{\overline{u}}
	=&-\aupw{-(u^m+1)^{n_2-1}\nabla\Pi_1 \log(u^m+\varepsilon)}{\widehat{u}_\oplus}{\overline{u}}\nonumber\\ &-\aupw{(u^m+1)^{n_2-1}\nabla v^{m+1}}{\widehat{u}_\oplus}{\bu}\nonumber\\&-\aupw{-(u^m+1)^{n_3-1}\nabla w^{m+1}}{\widehat{u}_\oplus}{\bu}\nonumber\\ &+\lambda \escalarL{(u^m)^\rho}{\bu} -\mu \escalarL{\widehat{u}_\oplus(u^m)^{k-1}}{\bu}\nonumber\\&- c \escalarL{\widehat{u}_\oplus (u^m)^{\gamma-1}\vert\nabla\Pi_1\log(u^m+\varepsilon)\vert^\gamma}{\bu},
	\end{align}
	for every $\bu\in\Pd_0(\T_h)$.
	
	It is straightforward to check that there is a unique solution $u$ of \eqref{eq:existence_solution_truncated} so $T$ is well defined.
	
	Secondly, we will check that $T$ is continuous. Let $\{\widehat{u}_j\}_{j\in\N}\subset\Pd_0(\T_h)$ be a sequence such that $\lim_{j\to\infty}\widehat{u}_j=\widehat{u}$. Taking into account that all norms are equivalent in $\Pd_0(\T_h)$ since it is a finite-dimensional space, the convergence $\widehat u_j\to \widehat u$ is equivalent to the elementwise convergence $(\widehat u_j)_K\to \widehat u_K$ for every $K\in\T_h$ (this may be seen, for instance, by using the norm $\norma{\cdot}_{L^\infty(\Omega)}$). Taking limits when $j\to \infty$ in \eqref{eq:existence_solution_truncated} (with $\widehat{u}\coloneqq\widehat{u}_j$ and $u\coloneqq T(\widehat u_j)$), using the notion of elementwise convergence,
	we get that $$\lim_{j\to \infty} T(\widehat u_j)=T(\widehat u)=T\left(\lim_{j\to \infty}\widehat u_j\right),$$ hence $T$ is continuous. In addition, $T$ is compact since $\Pd_0(\T_h)$ has finite dimension.
	
	Finally, let us prove that the set $$B=\{u\in\Pd_0(\T_h)\colon u=\alpha T(u)\text{ for some } 0\le\alpha\le1\}$$ is bounded (uniformly with respect to $\alpha$). The case $\alpha=0$ is trivial so we will assume that $\alpha\in(0,1]$.
	
	If $u\in B$, then $\frac{1}{\alpha}u\in\Pd_0(\T_h)$ is the solution of \eqref{eq:existence_solution_truncated} with $\widehat u=u$. Therefore, taking $\overline u=1$, we obtain $$\int_\Omega u\le\alpha\left(\int_\Omega u^m+\Delta t^{m+1}\lambda\int_\Omega (u^m)^\rho\right),$$ and, as $u^m\ge 0$ and it can be proved that $u\ge 0$ using the same arguments than in Theorem~\ref{thm:positivity_DG_chemo_truncated}, we get that $$\norma{u}_{L^1(\Omega)}\le \int_\Omega u^m+\Delta t^{m+1}\lambda\int_\Omega(u^m)^\rho.$$
	
	Since $\Pd_0(\T_h)$ is a finite-dimensional space where all the norms are equivalent, we have proved that $B$ is bounded.
	
	Thus, using the Leray-Schauder fixed point theorem, there is a fixed point of \eqref{eq:existence_solution_truncated}, therefore, we can conclude that there is a solution $u^{m+1}$ of
	\eqref{scheme_DG_upw_chemo_truncated}.
\end{proof}

Since every nonnegative solution of \eqref{scheme_DG_upw_chemo_truncated} is a solution of \eqref{scheme_DG_upw_chemo}, the following result is straightforward.

\begin{corollary}
    \label{cor:positivity-mass_scheme}
	There is at least one solution $u^{m+1}$ of \eqref{scheme_DG_upw_chemo} satisfying $u^{m+1}\ge 0$ and
	$$
	\int_\Omega u^{m+1}\le \int_\Omega u^m+\Delta t^{m+1} \lambda\int_\Omega (u^m)^\rho.
	$$
\end{corollary}

Notice that obtaining a nonnegative solution of the linear scheme \eqref{scheme_DG_upw_chemo} can be enforced by solving the
truncated nonlinear scheme \eqref{scheme_DG_upw_chemo_truncated}. However, in practice, the linear scheme \eqref{scheme_DG_upw_chemo} has provided a nonnegative approximation of the cell density in every numerical experiment that we have carried out without explicitly enforcing the nonnegativity constraint $u^{m+1}\ge 0$, as shown in Section~\ref{sec:numer-experiments}.

\begin{remark}
  Showing uniqueness of solution of the linear scheme \eqref{scheme_DG_upw_chemo} would imply that, as observed in the numerical tests, its solution is nonnegative without explicitly imposing the nonnegativity restriction $u^{m+1}\ge 0$. However, 
  this is not straightforward and it might require using inverse and trace inequalities that would probably involve some kind of restriction on the time step and mesh size, and this is left to a future work.
\end{remark}

%--------------------------------------
\section{Numerical experiments}
%--------------------------------------
\label{sec:numer-experiments}

Since the numerical tests shown in \cite{li2025dissipative} were confined to the parabolic-elliptic-elliptic version ($\tau=0$) of the problem \eqref{problem:chemotaxis_local}, we will focus in this work on the fully-parabolic version of \eqref{problem:chemotaxis_local} ($\tau=1$) and on the model \eqref{problem:chemotaxis_nonlocal}.

In this sense, we provide similar numerical examples to those in \cite{li2025dissipative} which are not only in accordance to the results in Section~\ref{sec:results}, but which also show that suppression of the condition \eqref{condgamma}, seems to lead to only local solution (i.e., $T$ finite) that blows up at $T$.

Here, we define $f_1(u)=u^\alpha$, $f_2(u)=u^\beta$ and we assume that all the parameters are set to $1$, but $k=1.1$, unless otherwise specified. For the sake of clarity, we also indicate in the figures any different value of the parameters with respect to those already fixed. 

  \begin{remark}
    Notice that the functions $f_1(u)=u^\alpha$, $f_2(u)=u^\beta$ do not necessarily satisfy Hypotheses~\ref{hyp:reglocal} and \ref{hyp:disf} for every $\alpha,\beta>0$. However, these functions are more physically meaningful than other similar choices of the type $f_1(u)=(u+\eta)^\alpha$, $f_2(u)=(u+\eta)^\beta$ for $\eta>0$ and, ultimately, 
    we expect the behavior of the solution to be similar to the case $\eta>0$ very small.
  \end{remark}

It is important to emphasize that the mass of the positive solutions $u^{m+1}$ of the scheme \eqref{scheme_DG_upw_chemo} is bounded for every $m\ge 0$, which implies $u^{m+1}\in L^1(\Omega)$ (Corollary~\ref{cor:positivity-mass_scheme}). Then, since all the norms are equivalent in a finite-dimensional space, we have as a consequence that $u^{m+1}\in L^\infty(\Omega)$ and, therefore, we cannot expect an actual blow-up in the discrete case as it occurs in the continuous model. However, we  observe in the numerical tests how the mass accumulates in some elements of $\T_h$ leading to the formation of peaks. In this sense,
 since an accurate discrete solution idealizing such blow-up scenario will exhibit a mass accumulation in small regions of the domain, 
we establish that blow-up occurs when, after that, the norm \mbox{$\norma{\cdot}_{L^\infty(\Omega)}$}  of the approximated solutions stabilizes over time (as in \cite{li2025dissipative,acosta-soba_KS_2022,shen_unconditionally_2020,badia2023bound}).

In order to compute the numerical tests, we have used the Python interface of the open source library \texttt{FEniCSx}
\cite{AlnaesEtal2014,ScroggsEtal2022, BasixJoss}
and the open source libraries \texttt{PyVista} \cite{sullivan2019pyvista} and \texttt{Matplotlib} \cite{Hunter:2007} to generate the plots. The source code for our implementation is hosted on GitHub\footnote{\url{https://github.com/danielacos/Papers-src}}.

In practice, in each time step, we have first computed a $\Pc_1(\T_h)$ approximation of the chemical signals $v^{m+1}$ and $w^{m+1}$ following the ideas in Section~\ref{sec:esquema_completamente_discreto}. Then, we have used the linear scheme \eqref{scheme_DG_upw_chemo} without explicitly enforcing the nonnegativity constraint $u^{m+1}\ge 0$, obtaining a positive approximation of the cell density, $u^{m+1}$, for every case below without needing any constraint on the time step or the mesh size. This strategy has allowed us to carry out the computationally demanding three-dimensional tests shown in the following section with considerably less computational effort than it would have been required if we had used a nonlinear approximation such as \eqref{scheme_DG_upw_chemo_truncated}.

Notice that, for the sake of a better visualization of the cell density in the 3D Figures~\ref{fig:test_1_p1c_u}, \ref{fig:test_3_p1c_u}, \ref{fig:test_3_2} (first row) and \ref{fig:test_3_3} (first row), we have plotted $\Pih_1 u^{m+1}$, i.e. the $\Pc_1(\T_h)$-projection of $u^{m+1}$ using mass-lumping as defined in \eqref{eq:esquema_DG_Pih1}. This projection preserves the positivity of the variable $u^{m+1}$ but introduces some diffusion that makes the maximum of $\Pih_1 u^{m+1}$ smaller than the actual maximum of $u^{m+1}$, compare for instance Figure~\ref{fig:test_1_p1c_u} with Figure~\ref{fig:test_1_max-u_a} in the case $c=0$.

\subsection{Only attraction local model}
\label{sec:test_1}

\begin{figure}[htbp]
  \centering
  \begin{subfigure}[b]{0.45\textwidth}
  \centering
      \includegraphics[width=0.9\textwidth]{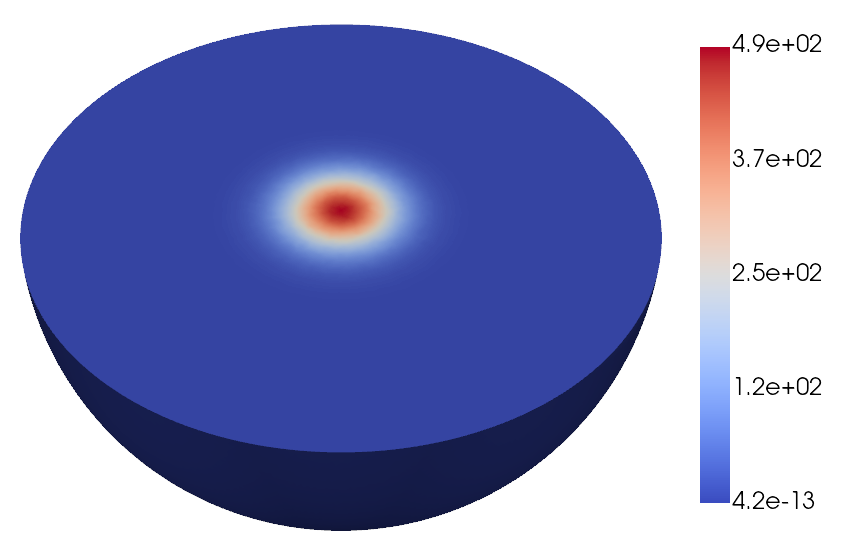}
      \caption{$t=0$}
      \label{fig:test_1_p1c_u_a}
  \end{subfigure}
  \hfill
  \begin{subfigure}[b]{0.45\textwidth}
  \centering
      \includegraphics[width=0.9\textwidth]{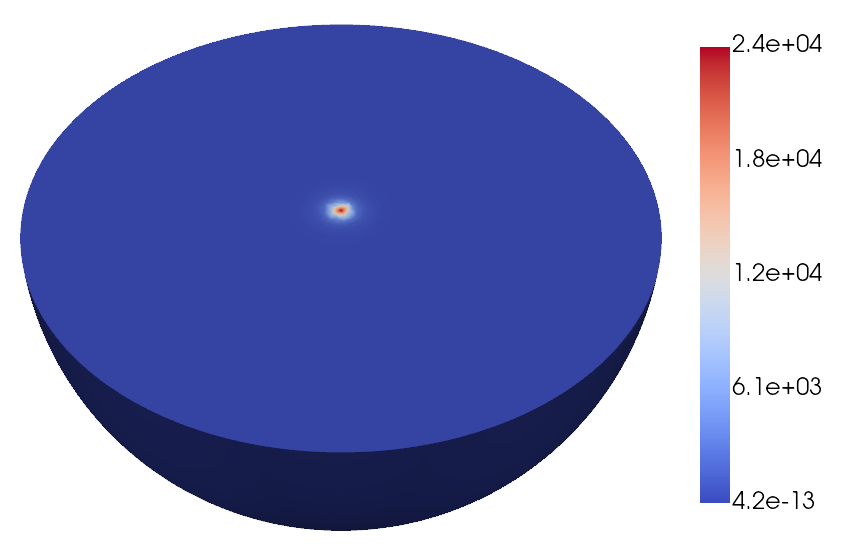}
      \caption{$t=5\cdot 10^{-4}$}
      \label{fig:test_1_p1c_u_b}
  \end{subfigure}
  \\
  \begin{subfigure}[b]{0.45\textwidth}
  \centering
      \includegraphics[width=0.9\textwidth]{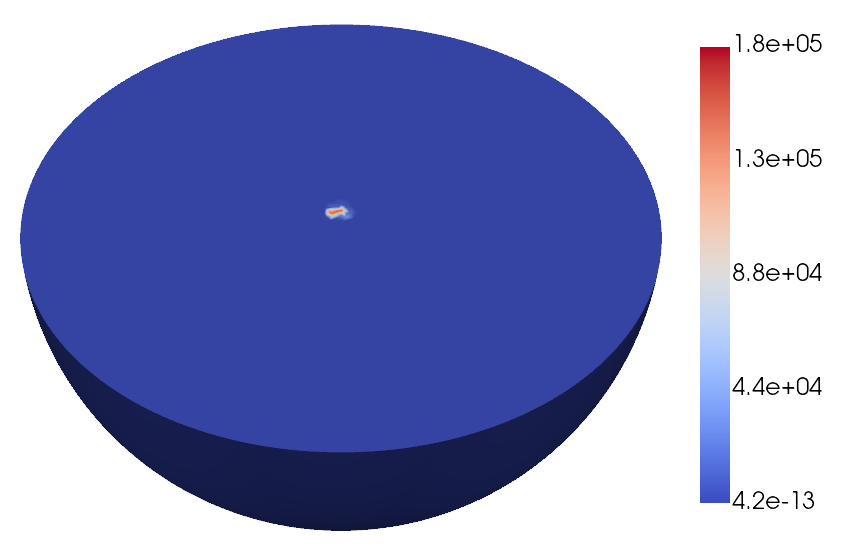}
      \caption{$t=10^{-3}$}
      \label{fig:test_1_p1c_u_c}
  \end{subfigure} 
  \hfill
  \begin{subfigure}[b]{0.45\textwidth}
  \centering
      \includegraphics[width=0.9\textwidth]{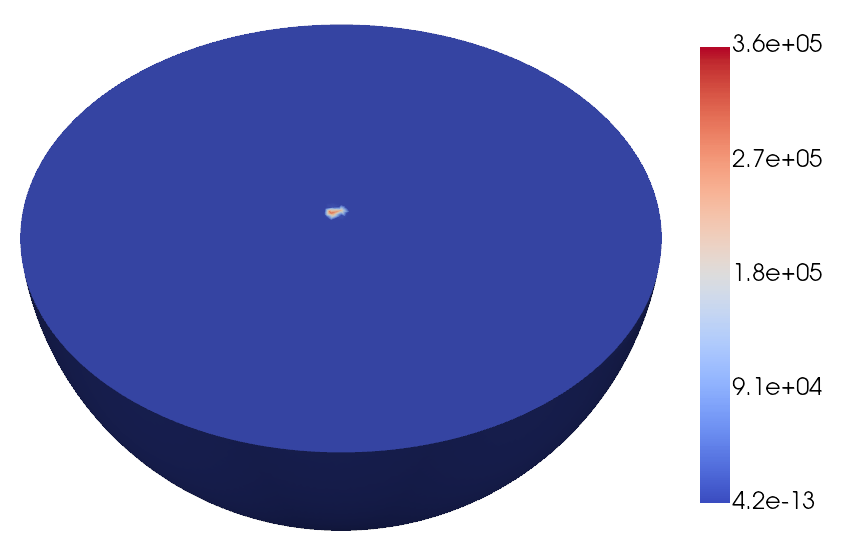}
      \caption{$t=3\cdot 10^{-3}$}
      \label{fig:test_1_p1c_u_d}
  \end{subfigure}      
  \caption{$\Pih_1 u^m$ at different time steps in Test~\ref{sec:test_1} ($c=0$, $\chi=5$, $\xi=0$)}
  \label{fig:test_1_p1c_u}
\end{figure}

\begin{figure}[htbp]
    \centering
    \begin{subfigure}[b]{0.45\textwidth}
    \centering
        \includegraphics[width=0.9\textwidth]{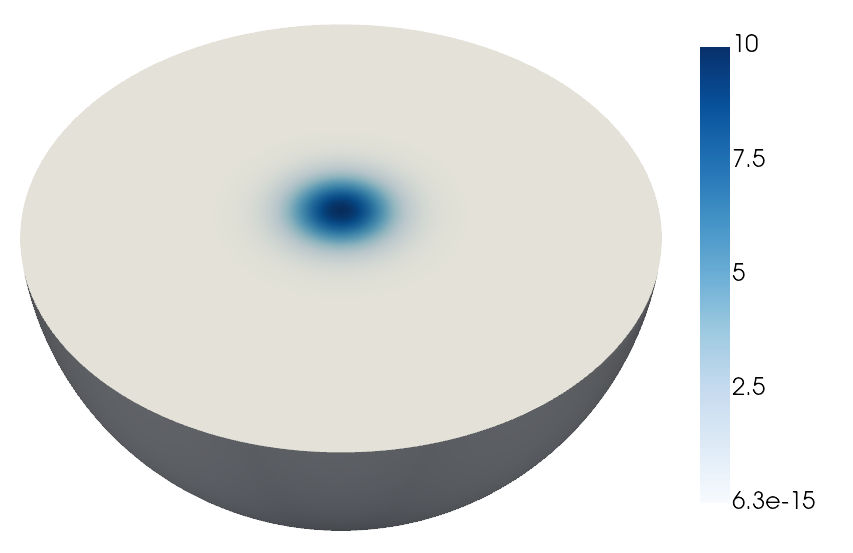}
        \caption{$t=0$}
        \label{fig:test_1_v_a}
    \end{subfigure}
    \hfill
    \begin{subfigure}[b]{0.45\textwidth}
    \centering
        \includegraphics[width=0.9\textwidth]{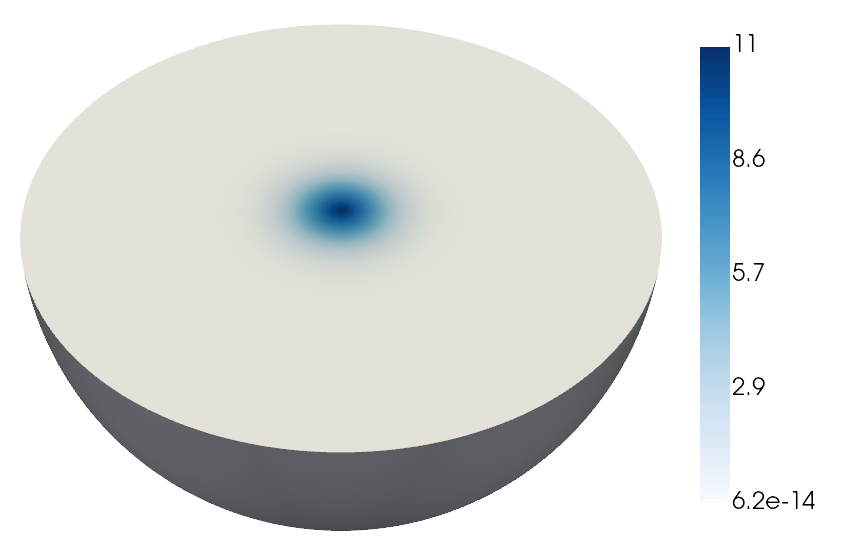}
        \caption{$t=5\cdot 10^{-4}$}
        \label{fig:test_1_v_b}
    \end{subfigure}
    \\
    \begin{subfigure}[b]{0.45\textwidth}
    \centering
        \includegraphics[width=0.9\textwidth]{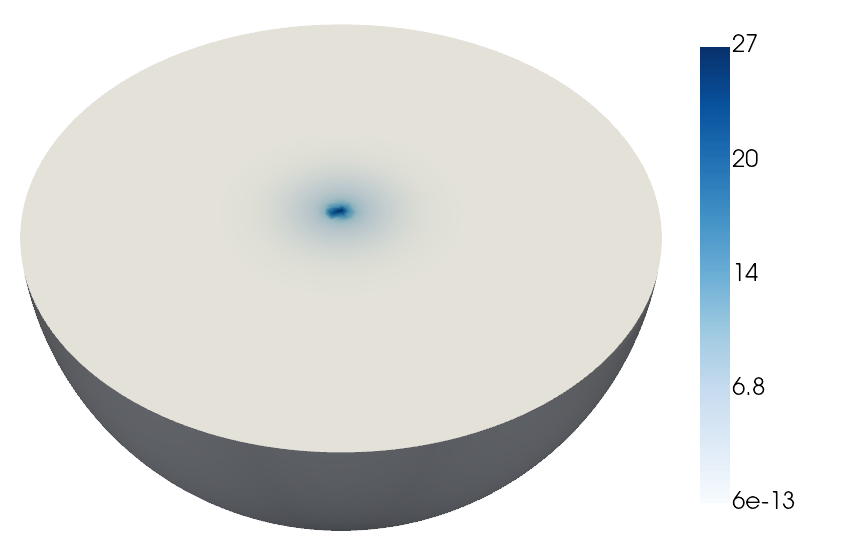}
        \caption{$t=10^{-3}$}
        \label{fig:test_1_v_c}
    \end{subfigure} 
    \hfill
    \begin{subfigure}[b]{0.45\textwidth}
    \centering
        \includegraphics[width=0.9\textwidth]{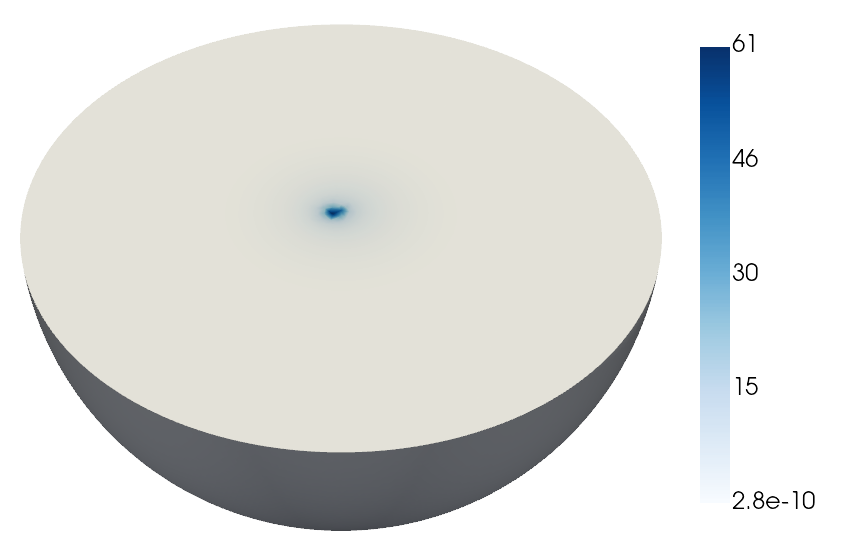}
        \caption{$t=3\cdot 10^{-3}$}
        \label{fig:test_1_v_d}
    \end{subfigure}      
    \caption{$v^{m}$ at different time steps in Test~\ref{sec:test_1} ($c=0$, $\chi=5$, $\xi=0$)}
    \label{fig:test_1_v}
  \end{figure}

In this first example, we are going to carry out a similar test to the first numerical experiment shown in \cite{li2025dissipative} for \eqref{problem:chemotaxis_local} with $\tau=0$, but, in this case, we will take $\tau=1$. For this test, we take a three-dimensional domain $\Omega=\{(x,y,z)\in\R^3\colon x^2+y^2+z^2< 1\}$ with the following initial condition used in \cite{li2025dissipative},
\begin{equation}
  \label{test_1_u_0}
  u_0(x,y,z)\coloneqq500e^{-35(x^2+y^2+z^2)},
\end{equation}
plotted in Figure~\ref{fig:test_1_p1c_u_a}. Notice that, in this section, for the sake of a better visualization of the results, we have only represented the half of the domain below the plane $z=0$. Moreover, we take $\chi=5$ and $\xi=0$ (no repulsion) and we use a mesh of size $h\approx4.4\cdot 10^{-2}$ and a constant time step $\Delta t^{m+1}=10^{-5}$ for every $m\ge 0$, which are the same parameters used in \cite{li2025dissipative}. Since $\tau=1$, $\chi=5$ and $\xi=0$ we need to add an initial condition for the chemoattractant density $v_0$, that we take as
\begin{equation}
  \label{test_1_v_0}
  v_0(x,y,z)\coloneqq10e^{-35(x^2+y^2+z^2)},
\end{equation}
which is plotted in Figure~\ref{fig:test_1_v_a}.

To the best of our knowledge, there are no theoretical results providing conditions for blow-up in the model \eqref{problem:chemotaxis_local} with $\tau=1$. Indeed, the existing literature on blow-up phenomena for fully parabolic chemotaxis models is very limited. While progress has been made in the simplified parabolic-elliptic cases, the techniques used there do not easily extend to fully parabolic models. In any case, we can observe in Figures~\ref{fig:test_1_p1c_u} and~\ref{fig:test_1_v} how the mass accumulates in the middle of the domain as in the test in \cite{li2025dissipative}.

\begin{figure}[htbp]
	\centering
	\begin{subfigure}{\textwidth}
	  \centering
	  \includegraphics[width=0.65\textwidth]{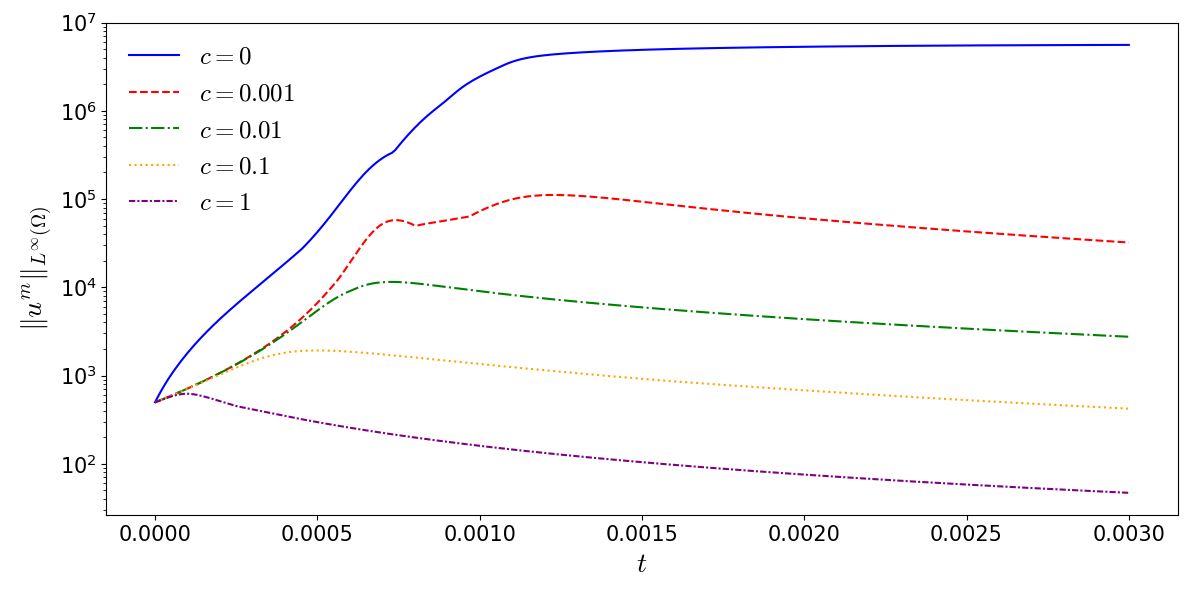}
	  \caption{$\gamma=1.75$}
	  \label{fig:test_1_max-u_a}
	\end{subfigure}
	\begin{subfigure}{\textwidth}
	  \centering
	  \includegraphics[width=0.65\textwidth]{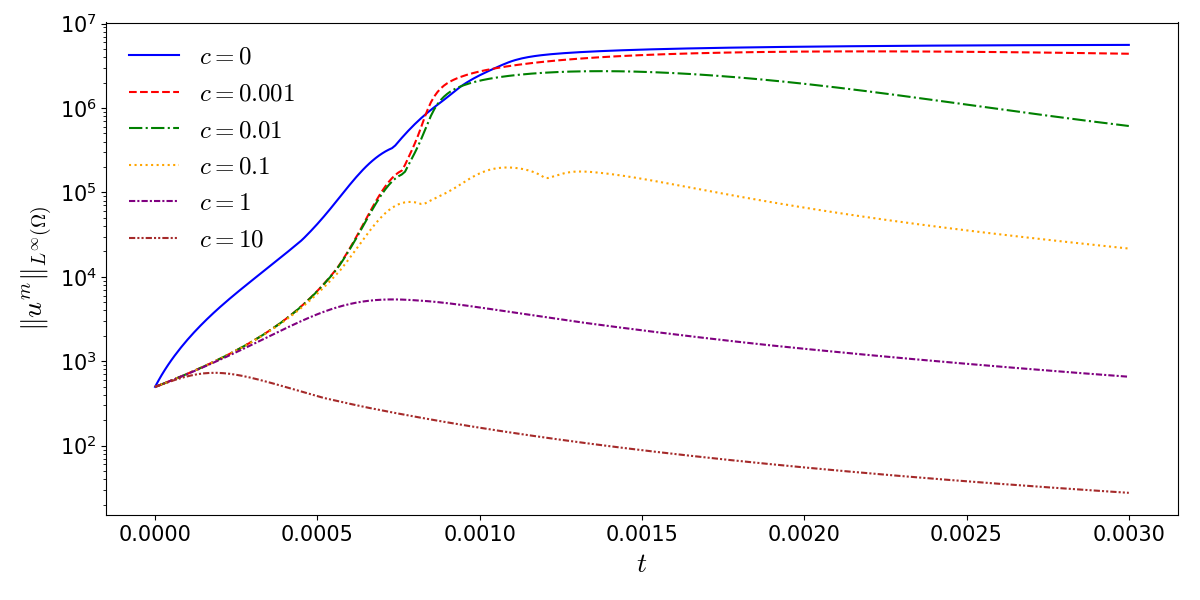}
	  \caption{$\gamma=1.4$}
	  \label{fig:test_1_max-u_b}
	\end{subfigure}
	\begin{subfigure}{\textwidth}
	  \centering
	  \includegraphics[width=0.65\textwidth]{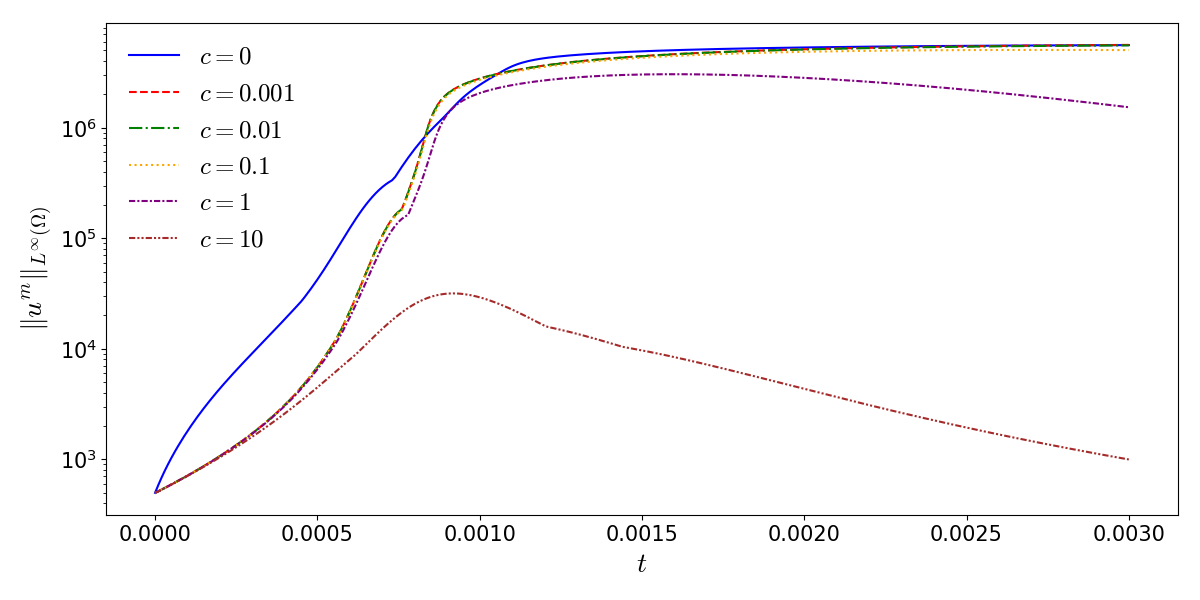}
	  \caption{$\gamma=1.1$}
	  \label{fig:test_1_max-u_c}
	\end{subfigure}
	\caption{$\|u ^m\|_{L^\infty(\Omega)}$ over time for different values of $c$ and $\gamma$ in Test~\ref{sec:test_1} ($\chi = 5$, $\xi = 0$)}
	\label{fig:test_1_max-u}
  \end{figure}

This mass-accumulation can be prevented using the damping gradient as shown in Theorem~\ref{thm:theoremlocal}. In this sense, we show in Figure~\ref{fig:test_1_max-u} the evolution of the maximum of $u^m$ for different choices of $c$ and $\gamma$. In this case, where $\chi=0$, $n_1=\alpha=1$, the condition \eqref{condgamma} can be reduced to $\frac{2d}{d+1}\le \gamma \le 2$, see \cite{ishida2023keller}. Hence, as expected, this mass accumulation is prevented for whatever choice of $c>0$ (even for small values like $c=10^{-3}$) if $1.5<\gamma\le 2$ (see Figure~\ref{fig:test_1_max-u_a}). For values $1\le\gamma\le1.5$, a sufficiently big value $c>0$ is needed to prevent this mass accumulation, see Figures~\ref{fig:test_1_max-u_b} and~\ref{fig:test_1_max-u_c}. Notice that these results are similar to those shown in \cite{li2025dissipative} for $\tau=0$.

\subsection{Attraction-repulsion local model}
\label{sec:test_2}
Now, a similar experiment to the previous numerical test is carried out for the model \eqref{problem:chemotaxis_local} where, again $\tau=1$ and $\chi=5$ but, in this case, we consider the chemorepulsion and take $\xi=1$. We keep the same initial conditions $u_0$ and $v_0$ than in the previous test, \eqref{test_1_u_0} and \eqref{test_1_v_0}, respectively, and the discrete parameters $h\approx4.4\cdot 10^{-2}$, $\Delta t^{m+1}=10^{-5}$ for every $m\ge 0$. However, in this case, we need to add an initial condition for the chemorepulsive density $w_0$ that we define as
\begin{equation}
  \label{test_2_w_0}
  w_0(x,y,z)\coloneqq10e^{-35(x^2+y^2+z^2)},
\end{equation}
so that $w_0=v_0$ in $\Omega$.

\begin{figure}[H]
  \centering
  \begin{subfigure}{\textwidth}
    \centering
    \includegraphics[width=0.65\textwidth]{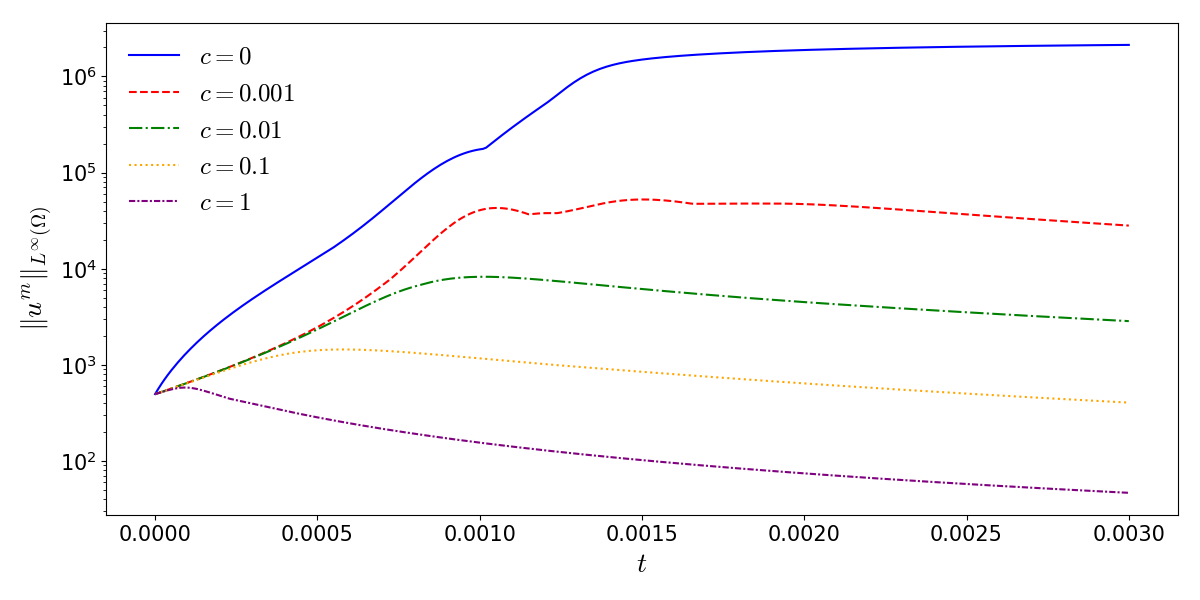}
    \caption{$\gamma=1.75$}
    \label{fig:test_2_max-u_a}
  \end{subfigure}
  \begin{subfigure}{\textwidth}
    \centering
    \includegraphics[width=0.65\textwidth]{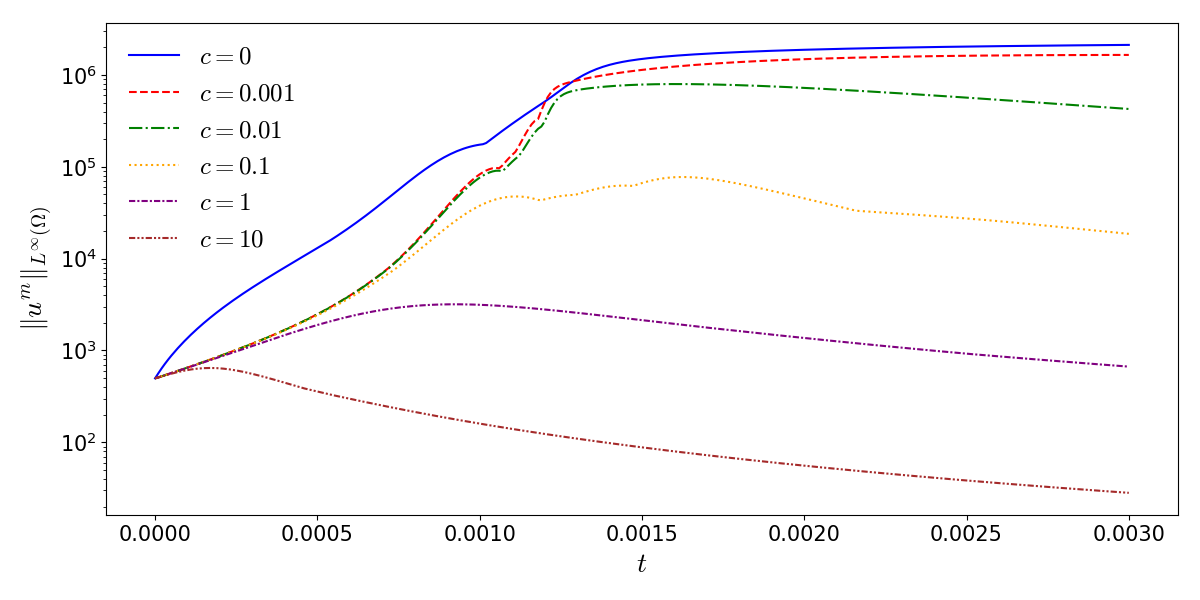}
    \caption{$\gamma=1.4$}
    \label{fig:test_2_max-u_b}
  \end{subfigure}
  \begin{subfigure}{\textwidth}
    \centering
    \includegraphics[width=0.65\textwidth]{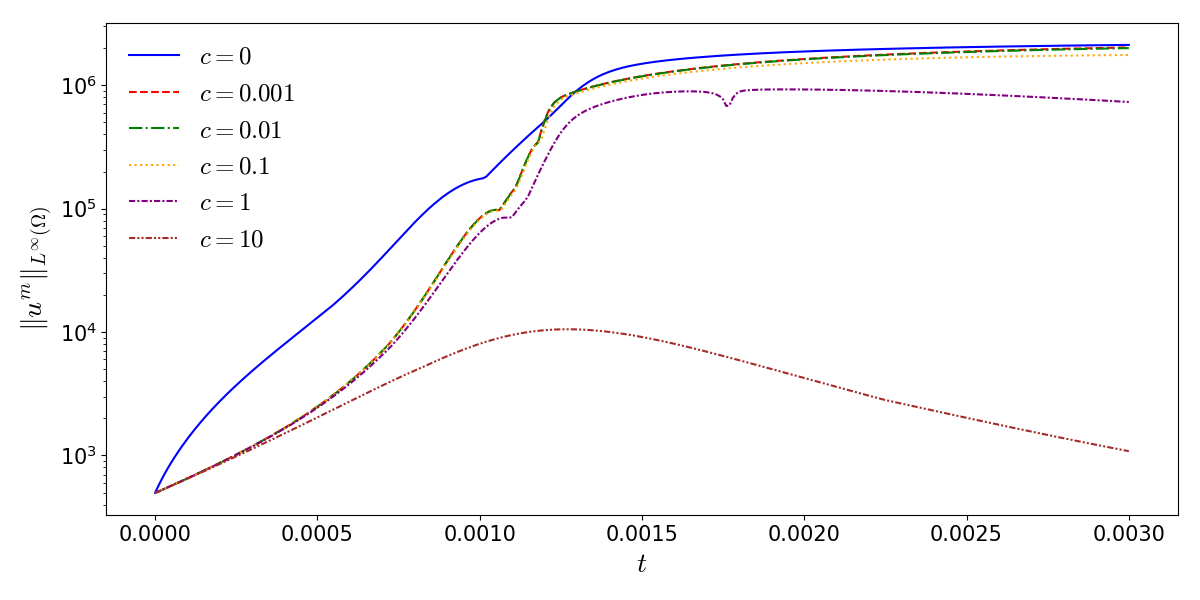}
    \caption{$\gamma=1.1$}
    \label{fig:test_2_max-u_c}
  \end{subfigure}
  \caption{$\|u ^m\|_{L^\infty(\Omega)}$ over time for different values of $c$ and $\gamma$ in Test~\ref{sec:test_2} ($\chi = 5$, $\xi = 1$)}
  \label{fig:test_2_max-u}
\end{figure}

In Figure~\ref{fig:test_2_max-u} we have plotted the evolution of the maximum of $u^m$. If we compare these results with those for the previous test in Figure~\ref{fig:test_1_max-u}, we observe a similar behavior, where a mass accumulation appears when $c=0$. This mass accumulation is prevented for $\gamma=1.75$ independently of $c>0$ as it satisfies \eqref{condgamma}, i.e. $1.5\le\gamma\le 2$, while for $\gamma=1.4$ and $\gamma=1.1$, a big enough value of $c>0$ is needed. However, in this case, the mass accumulation seems to occur later and be less pronounced than in the previous test due to the effect of the chemorepulsion.

\subsection{Attraction-repulsion nonlocal model}
\label{sec:test_3}

\begin{figure}[htbp]
  \centering
  \begin{subfigure}[b]{0.45\textwidth}
  \centering
      \includegraphics[width=0.9\textwidth]{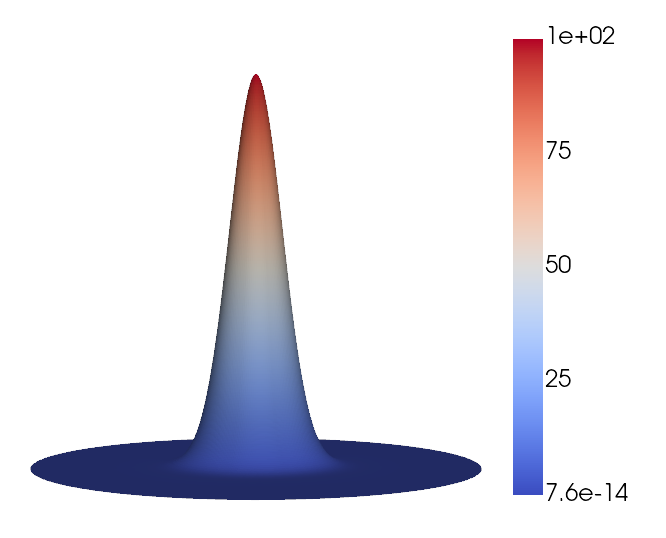}
      \caption{$t=0$}
      \label{fig:test_3_p1c_u_a}
  \end{subfigure}
  \hfill
  \begin{subfigure}[b]{0.45\textwidth}
  \centering
      \includegraphics[width=0.9\textwidth]{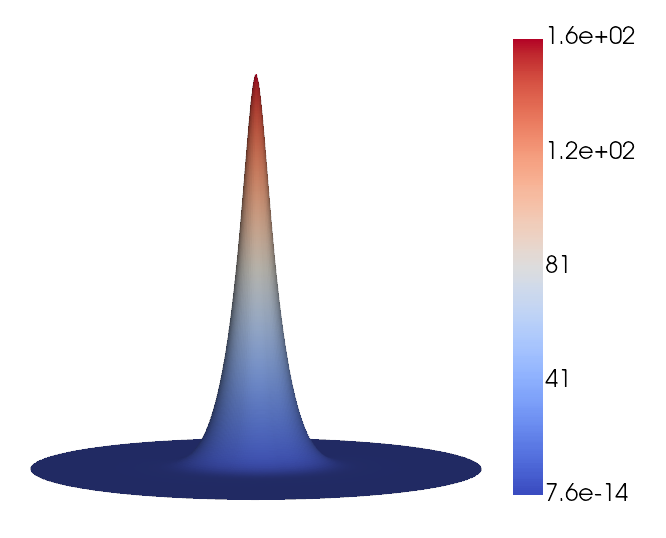}
      \caption{$t=5\cdot 10^{-4}$}
      \label{fig:test_3_p1c_u_b}
  \end{subfigure}
  \\
  \begin{subfigure}[b]{0.45\textwidth}
  \centering
      \includegraphics[width=0.9\textwidth]{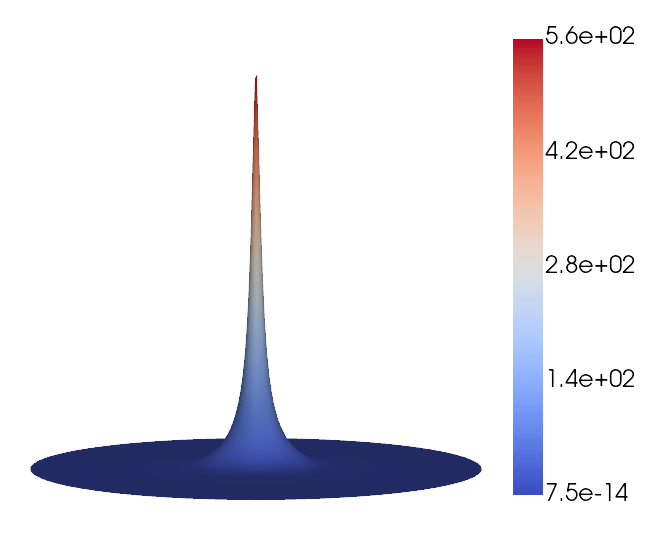}
      \caption{$t=10^{-3}$}
      \label{fig:test_3_p1c_u_c}
  \end{subfigure} 
  \hfill
  \begin{subfigure}[b]{0.45\textwidth}
  \centering
      \includegraphics[width=0.9\textwidth]{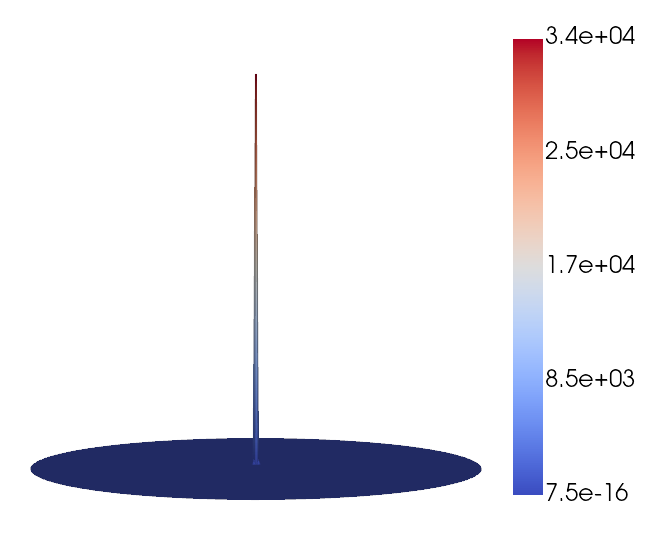}
      \caption{$t=3\cdot 10^{-3}$}
      \label{fig:test_3_p1c_u_d}
  \end{subfigure}      
  \caption{$\Pih_1 u^{m}$ at different time steps in Test~\ref{sec:test_3} ($c=0$, $\alpha=1.5$)}
  \label{fig:test_3_p1c_u}
\end{figure}

Now, we consider the nonlocal model \eqref{problem:chemotaxis_nonlocal} and compute a numerical experiment in the two-dimensional domain $\Omega=\{(x,y)\in\R^2\colon x^2+y^2<1\}$. We take the initial condition for the cell density
$$
u_0(x,y,z)\coloneqq100e^{-35(x^2+y^2)},
$$
plotted in Figure~\ref{fig:test_3_p1c_u_a}, and set $\alpha=1.5$ and the discrete parameters $\Delta t^{m+1}=10^{-5}$ for every $m\ge 0$ and $h\approx1.4\cdot 10^{-2}$. Moreover, since $\chi,\xi>0$, we need to add an initial condition for the chemical signals $v_0$ and $w_0$ that we take as
$$
  v_0(x,y,z)\coloneqq10e^{-35(x^2+y^2)}, \quad w_0(x,y,z)\coloneqq10e^{-35(x^2+y^2)},
$$
so that $v_0=w_0$ in $\Omega$.

With this choice, we are under the assumptions in \cite[Theorem 3]{ColumbuFrassuViglialoro2023}, where for $\chi,\xi,\lambda,\mu, n_2=n_3>0$, $\rho=1$, $k>1$ and
$$
\alpha>\beta \ \mathrm{ and }\ 
\begin{cases}
n_2+\alpha>\max\{n_1+\frac{d}{2}k,k\} & \mathrm{if\ } n_1\ge 0,\\
n_2+\alpha>\max\{\frac{d}{2}k,k\} & \mathrm{if\ } n_1< 0,
\end{cases}
$$
that guarantee the existence of a blowing-up solution in the case $c=0$ for a certain initial condition that satisfies
$
\int_\Omega u_0\ge \left(\frac{\lambda}{\mu}\vert\Omega\vert^{k-1}\right)^{\frac{1}{k-1}}.
$
In particular, the numerical approximation of the solution, plotted in Figure~\ref{fig:test_3_p1c_u}, exhibits a mass accumulation in the middle of the domain that suggests that a chemotactic collapse has occurred.

\begin{figure}[htbp]
\centering

\begin{subfigure}[b]{0.45\textwidth}
	\centering
	\hspace*{-2cm}
	(a) $\gamma=1.4$
\end{subfigure}
\hfill
% \hspace*{-1cm}
\begin{subfigure}[b]{0.45\textwidth}
	\centering
	\hspace*{-2cm}
	(b) $\gamma=1.75$
\end{subfigure}

\begin{subfigure}[b]{0.45\textwidth}
	\centering
	\includegraphics[width=0.9\textwidth]{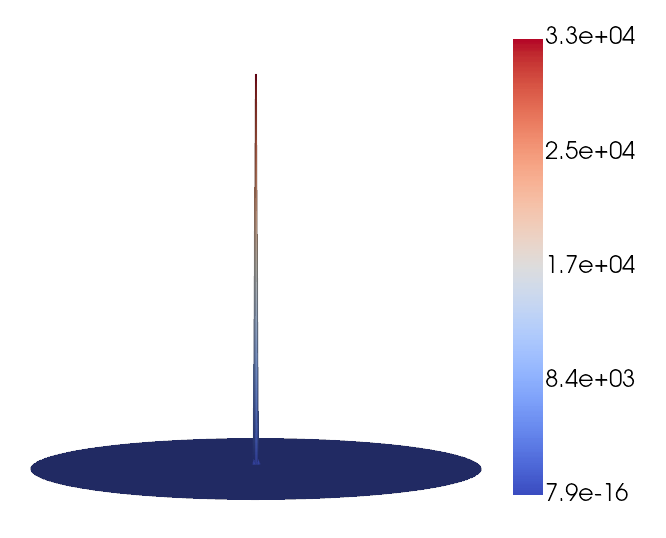}
\end{subfigure}\hfill
\begin{subfigure}[b]{0.45\textwidth}
	\centering
	\includegraphics[width=0.9\textwidth]{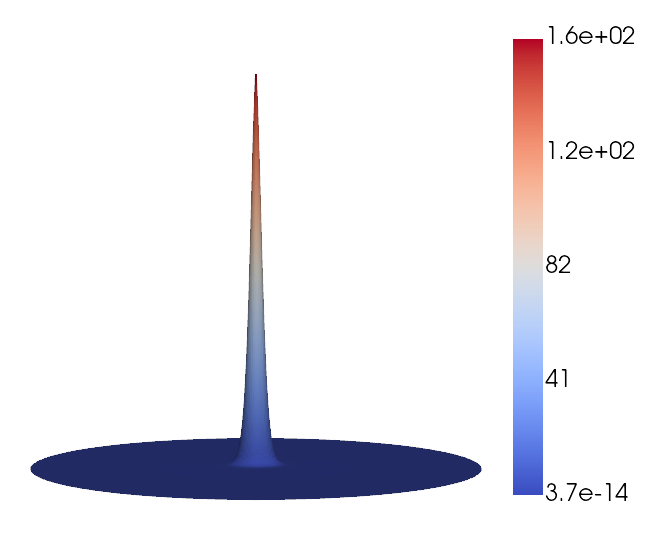}
\end{subfigure}

\begin{subfigure}[b]{0.45\textwidth}
	\centering
	\includegraphics[width=0.9\textwidth]{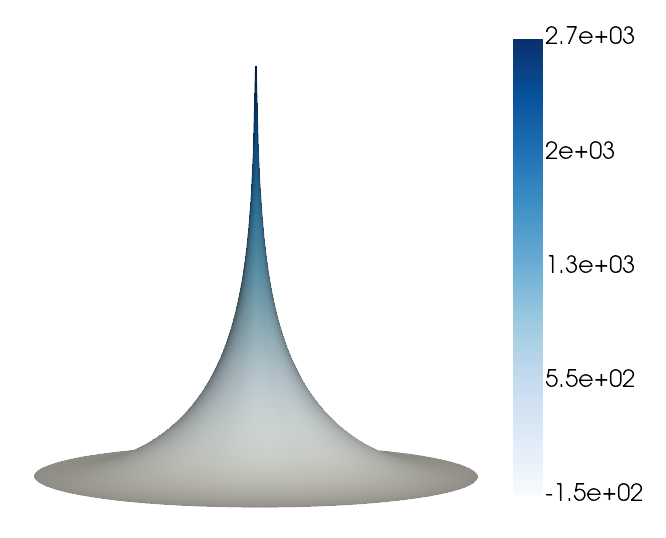}
\end{subfigure}\hfill
\begin{subfigure}[b]{0.45\textwidth}
	\centering
	\includegraphics[width=0.9\textwidth]{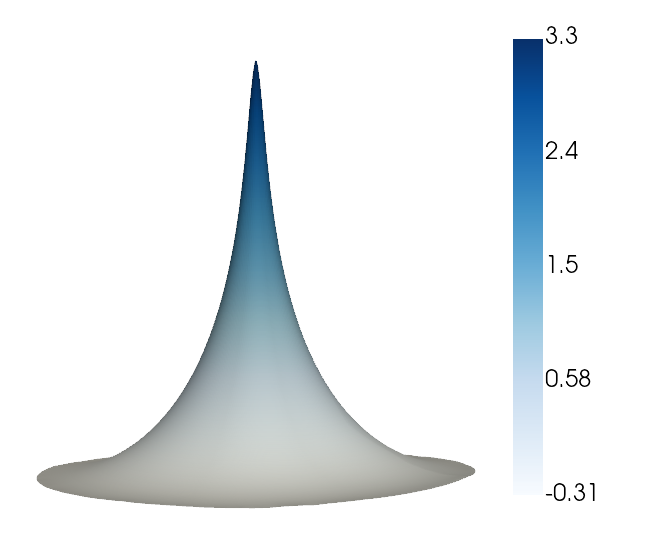}
\end{subfigure}

\begin{subfigure}[b]{0.45\textwidth}
	\centering
	\includegraphics[width=0.9\textwidth]{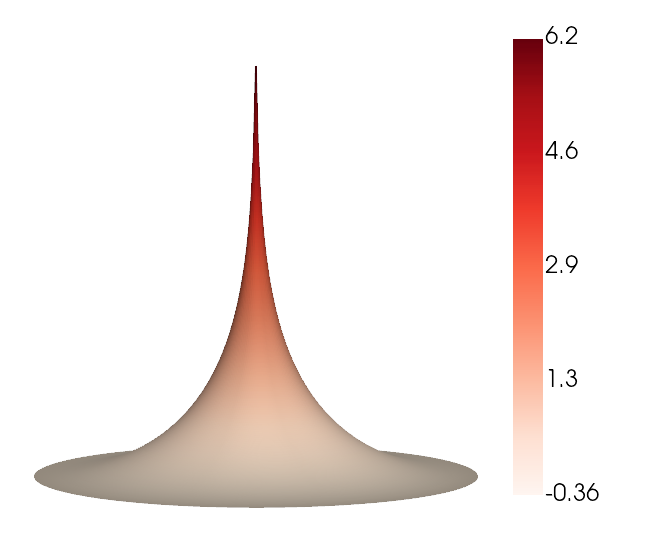}
\end{subfigure}\hfill
\begin{subfigure}[b]{0.45\textwidth}
	\centering
	\includegraphics[width=0.9\textwidth]{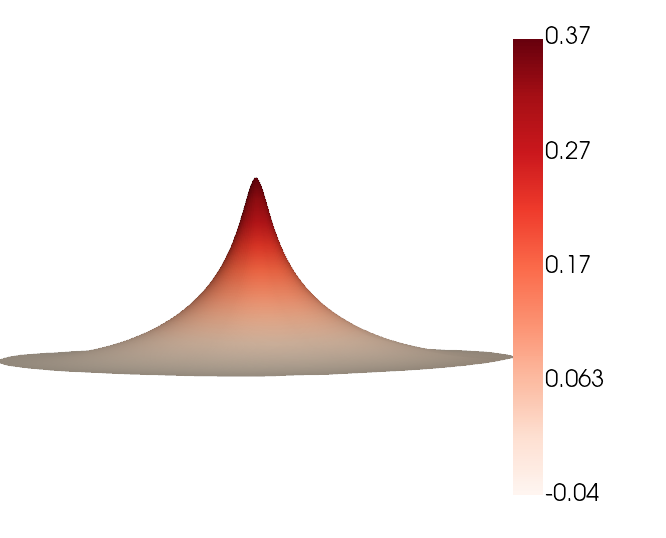}
\end{subfigure}
\caption{$\Pih_1 u^{m}$ (first row), $v^{m}$ (second row) and $w^{m}$ (third row) at $t=3\cdot 10^{-3}$ for different values of $\gamma$ and $c=10^{-3}$ in Test~\ref{sec:test_3} ($\alpha=1.5$)}
\label{fig:test_3_2}
\end{figure}

\begin{figure}[htbp]
    \centering
    \begin{subfigure}{0.45\textwidth}
		\centering
		\hspace*{-2cm}
        (a) $\gamma=1.4$
	\end{subfigure}\hfill
	\begin{subfigure}{0.45\textwidth}
		\centering
		\hspace*{-2cm}
		(b) $\gamma=1.75$
	\end{subfigure}

	\begin{subfigure}{0.45\textwidth}
		\centering
        \includegraphics[width=0.9\textwidth]{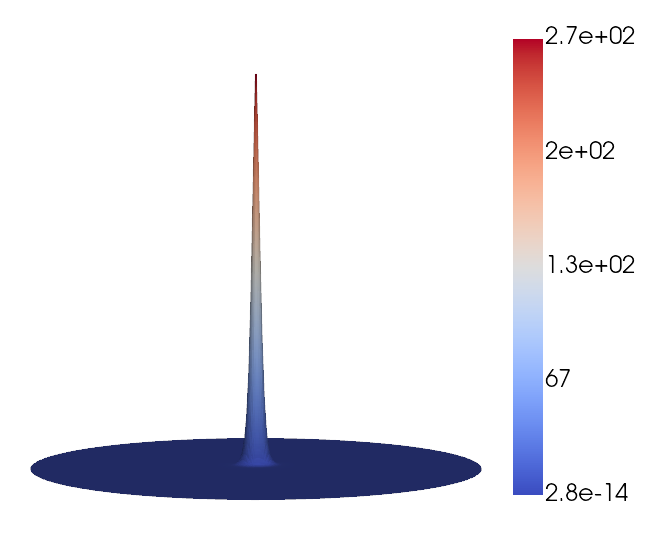}
	\end{subfigure}\hfill
	\begin{subfigure}{0.45\textwidth}
		\centering
        \includegraphics[width=0.9\textwidth]{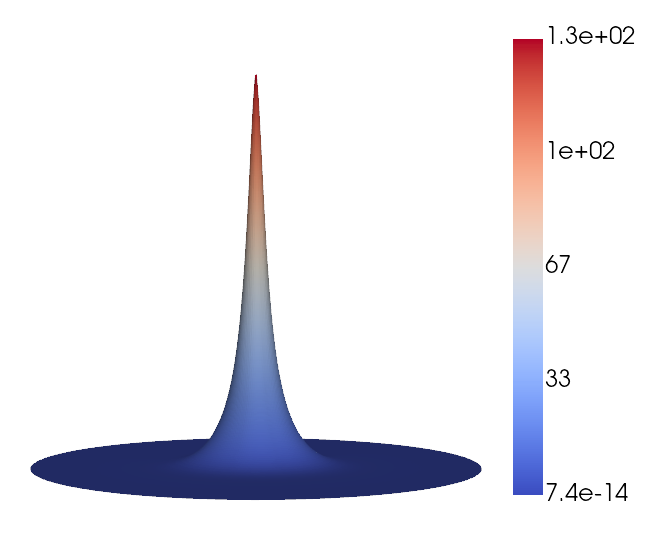}
	\end{subfigure}

	\begin{subfigure}{0.45\textwidth}
		\centering
        \includegraphics[width=0.9\textwidth]{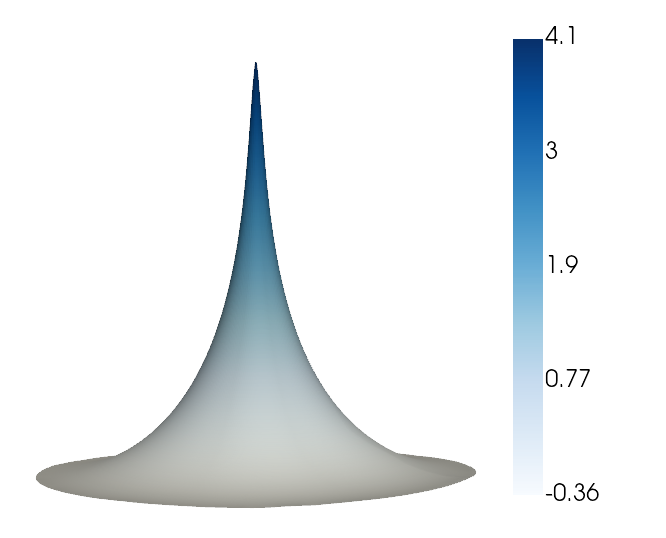}
	\end{subfigure}\hfill
	\begin{subfigure}{0.45\textwidth}
		\centering
        \includegraphics[width=0.9\textwidth]{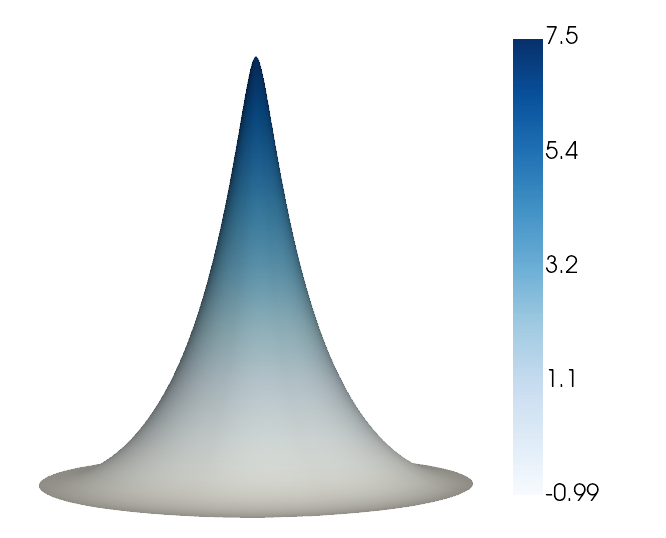}
	\end{subfigure}

	\begin{subfigure}{0.45\textwidth}
		\centering
        \includegraphics[width=0.9\textwidth]{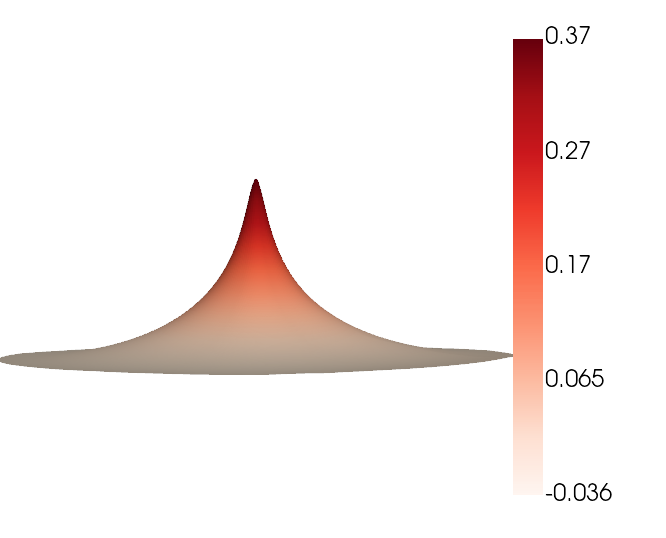}
	\end{subfigure}\hfill
	\begin{subfigure}{0.45\textwidth}
		\centering
        \includegraphics[width=0.9\textwidth]{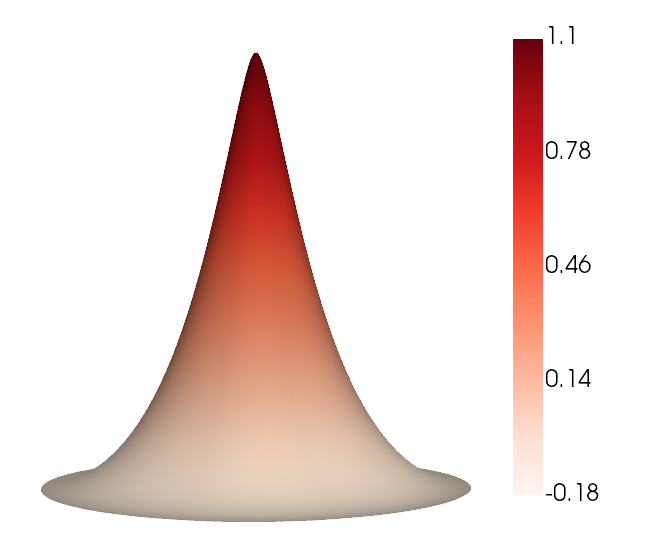}
    \end{subfigure}
    \caption{$\Pih_1 u^{m}$ (first row), $v^{m}$ (second row) and $w^{m}$ (third row) at $t=3\cdot 10^{-3}$ for different values of $\gamma$ and $c=10^{-1}$ in Test~\ref{sec:test_3} ($\alpha=1.5$)}
    \label{fig:test_3_3}
\end{figure}

According to Theorem~\ref{thm:theoremnonlocal}, if \eqref{condgamma} is satisfied, i.e., $1.67\approx5/3<\gamma\le 2$, then the blow-up is prevented for any value $c>0$. In fact, a comparison of the different results obtained for different values of $\gamma\in[1,2]$ and $c>0$ are shown in Figures~\ref{fig:test_3_2} and~\ref{fig:test_3_3}. Notice how the mass accumulation does no longer appear for the value $\gamma=1.75$, which satisfies \eqref{condgamma}. However, for $\gamma=1.4< 5/3\approx 1.67$, a big enough value of $c>0$ seems to be required to avoid the singularity.

\section*{Acknowledgments}
  The first author has been supported by UCA FPU contract UCA/REC14VPCT/ 2020 funded by Universidad de Cádiz, by mobility grants funded by Plan Propio - UCA 2022-2023 and by a Graduate Scholarship funded by The University of Tennessee at Chattanooga.
  
  The second author is a member of the {\em Gruppo Nazionale per l'Analisi Matematica, la Probabilità e le loro Applicazioni} (GNAMPA) of the Istituto Nazionale di Alta Matematica (INdAM). This author is also supported by GNAMPA-INdAM Project \textit{Problemi non lineari di tipo stazionario ed evolutivo} (CUP--E53C23001670001).
  
  The third author has been supported by Grant PGC2018-098308-B-I00 (MCI/AEI/FEDER, UE, Spain), Grant US-1381261 (US/JUNTA/FEDER, UE, Spain) and Grant P20-01120 (PAIDI/JUNTA/ FEDER, UE, Spain).

\printbibliography

\end{document}